\documentclass[12pt, leqno]{amsart}

\usepackage[OT2,T1]{fontenc}
\usepackage{indentfirst}
\usepackage{amstext}
\usepackage{amsthm}
\usepackage{amsopn}
\usepackage{amsfonts}
\usepackage{amsmath}
\usepackage{latexsym}
\usepackage[francais,english]{babel}
\usepackage{amscd}
\usepackage{amssymb}
\usepackage{amsmath}
\usepackage[all,cmtip]{xy}
\usepackage{leftidx}
\usepackage{graphicx}
\usepackage{etoolbox}
\patchcmd{\section}{\normalfont\scshape\centering}{\normalfont\bfseries}{}{}
\patchcmd{\subsection}{-.5em}{.5em}{}{}

\newtheorem{theo}{{Theorem}}[section]
\newtheorem{coro}[theo]{{Corollary}}
\newtheorem{lemma}[theo]{{Lemma}}
\newtheorem{prop}[theo]{Proposition}

\theoremstyle{definition}
\newtheorem{remark}[theo]{\textbf{Remark}}
\newtheorem{defn}[theo]{Definition}
\newtheorem{example}[theo]{Example}
\numberwithin{equation}{section}

\newtheorem{notation}[theo]{Notation}

\DeclareFontEncoding{OT2}{}{} 

\begin{document}
\tolerance 400 \pretolerance 200 \selectlanguage{english}

\title{$K3$ surfaces with maximal complex multiplication}
\author{Eva  Bayer-Fluckiger}
\date{\today}
\maketitle

\begin{abstract} Let $X$ be a complex projective $K3$ surface having complex multiplication by a CM field $E$, and let $T_X$ be its
transcendental lattice. We say
that $X$ has {\it maximal complex multiplication} if ${\rm End}_{\rm Hdg}(T_X)$ is the ring of integers of $E$. 

\medskip For which CM fields $E$ does such a $K3$ surface exist ? What are the possibities for the transcendental
lattices, Picard lattices of these surfaces~?  The aim of this paper is to study these questions and give some examples.




\medskip

\end{abstract}

\small{} \normalsize

\medskip

\selectlanguage{english}
\section{Introduction}

Let $E$ be a CM field, and let $O_E$ be its ring of integers. A complex projective $K3$ surface $X$ is said to have  {\it complex multiplication by $E$} 
(or CM by $E$, for short) if 
${\rm End}_{\rm Hdg}(T_X \otimes_{\bf Z}{\bf Q}) \simeq E$ and ${\rm rank}(T_X) = [E: \bf {Q}]$, where $T_X$ is the transcendental
lattice of $X$. This implies that $[E:{\bf Q}] \leqslant 20$, and Taelman proved that if $[E:{\bf Q}] \leqslant 20$ then
there exists a $K3$ surface with CM by $E$ (cf. \cite {T}, Theorem 3).

\medskip Following Valloni \cite{V 21}, we say that a $K3$ surface $X$ has {\it complex multiplication by $O_E$} if $X$ has CM by $E$ and moreover
${\rm End}_{\rm Hdg}(T_X) \simeq O_E$; we also say that $X$ has then {\it maximal complex multiplication} (in \cite{V 21} this is called  ``principal complex multiplication'').

\medskip
If $X$ is a $K3$ surface, set $L_X = H^2(X,{\bf Z})$; the intersection form makes $L_X$ into a (unimodular) lattice, and $T_X$ is a sublattice of $L_X$. 
If moreover
$X$ has maximal complex multiplication, then $T_X$ has a structure of $O_E$-module, hence so has its dual $T_X^{\sharp}$; 
therefore the quotient $T_X^{\sharp}/T_X$ is isomorphic to $O_E/{\mathcal D}_X$, where $\mathcal D_X \subset O_E$ is an $O_E$-ideal,
called the {\it discriminant ideal} of $X$; note that the norm of $\mathcal D_X$ is the determinant of $T_X$, hence also the
absolute value of the determinant of the Picard lattice of $X$. 

\medskip
{\bf Question 1.} What are the possibilities for the ideal $\mathcal D_X$ ? 

\medskip
This is the subject matter of \S \ref{given section} (see Corollary \ref{D existence}), based on results of \S \ref{discriminant form} and \S \ref{K3 section}.

\medskip
The next issue is to classify up to isomorphism the $K3$ surfaces
with a given discriminant ideal : 

\medskip
{\bf Question 2.} What are the possibilities for the $K3$ surfaces $X$ with a given discriminant ideal $\mathcal D_X$ ?

\medskip
It is well-known that if $E$ is an imaginary quadratic field, then the isomorphism classes of elliptic curves with complex multiplication
by $O_E$ are in bijection with the ideal class group of $E$. 

\medskip
For an arbitrary CM field $E$, 
we define a finite group $\mathcal C$ having a similar property for the isomorphism classes of $K3$ surfaces with complex
multiplication by $O_E$ with the same discriminant ideal  (see Corollary \ref{K3 with D}).  We show

\medskip
\noindent
{\bf Theorem 1.} {\it There are only finitely many isomorphism classes of $K3$ surfaces having  maximal complex multiplication by the
same CM field and the same discriminant ideal.}

\medskip

The next sections contain some applications of the previous results.  In \cite {V 21}, Valloni raised the following question :

\medskip
{\bf Question 3.} For which CM fields $E$ do there exist $K3$ surfaces with maximal CM by $E$ ?

\medskip
 Valloni proved that if $[E:{\bf Q}] \leqslant 10$, then there exist infinitely many non-isomorphic $K3$ surfaces having CM by $O_E$.
The aim of \S \ref{CM section} is to give a sufficient criterion for the existence of infinitely many non-isomorphic $K3$ surfaces with CM by $O_E$ in terms of ramification properties of $E$  (see Theorem \ref{CM theorem}). One of the applications is a generalization of Valloni's result (cf. Corollary \ref{14}) :

\medskip
\noindent
{\bf Theorem 2.} {\it  If $[E:{\bf Q}] \leqslant 14$, then there exist infinitely many non-isomorphic complex projective  $K3$ surfaces with
complex multiplication by $O_E$.}

\medskip
This is no longer true in general if $[E:{\bf Q}] = 16, 18$ or $20$
(cf. Proposition \ref{high ramification} and Example \ref{counterexample}); but it does hold for {\it cyclotomic fields} :

\medskip
\noindent
{\bf Theorem 3.} {\it If $E$ is a cyclotomic field with  $2 \leqslant [E:{\bf Q}] \leqslant 20$, then there exist infinitely many non-isomorphic complex projective  $K3$ surfaces with complex multiplication by $O_E$.}

\medskip
In another direction, maximal complex multiplication on the transcendental lattice implies some properties
of the Picard lattice, such as possible degrees of polarisation, existence of elliptic fibrations. 

\medskip
A start on this is made in \S \ref{Picard section}; the
following example illustrates the results of this section (see Example \ref{44 and 66}; here $U$ denotes the 2-dimensional hyperbolic lattice, and for all integers $N$, the
lattice $U(N)$ denotes $U$ multiplied by $N$) :

\medskip
{\bf Example 1.} Let $E = {\bf Q}(\zeta_m)$ with $m = 44$ or $66$. There exists a $K3$ surface with maximal complex multiplication by $E$ 
with Picard lattice $L$ $\iff$ $L \simeq U(N)$ where $N \geqslant 1$ is an integer  $\equiv 1 \ {\rm (mod \ {\it m})}$
such that  all the prime divisors of $N$ are  $ \equiv \pm 1 \ {\rm (mod \ {\it m})}$.

\bigskip With this strategy in mind,
the second part of the paper concerns $K3$ surfaces having maximal complex multiplication by cyclotomic fields. One of
the results is the following

\medskip
\noindent
{\bf Theorem 4.} {\it Let $p$ be an odd prime number with $3 \leqslant p \leqslant 11$, and set $E = {\bf Q}(\zeta_p)$. Let $a \geqslant 1$ be an  odd integer.
 There exists a unique {\rm (up to isomorphism)} complex projective
 $K3$ surface $X_a$ with maximal complex multiplication by $E$ such that 
 ${\rm det}(T_{X_a}) = p^a$. 
 
 \medskip Moreover, the surfaces $X_a$ are isogeneous for all $a \geqslant 1$.}
 
 \medskip If $a = 1$, these surfaces  are isomorphic to Vorontsov's $K3$ surfaces (see \cite{V}, \cite {Ko}).
 
 \medskip
 If $p = 3, 7$ or $11$,  then for all $a \geqslant 1$ the $K3$ surfaces $X_a$ above and their twists (cf.  Definition \ref{K3 twist}) have automorphisms of order $p$ inducing the complex multiplication by $O_E$ (see Theorem \ref{twisted K3 lemma}),
 and the corresponding pairs ($K3$ surface, automorphism of order $p$) are the CM points by $O_E$ of the moduli spaces ${\mathcal M}^p_{K3}$ defined
 by Artebani, Sarti and Taki in \cite{AST} (see \S \ref{moduli}); see also \cite{AS} for $p = 3$ and \cite{OZ}, \cite{ACV} for $p = 11$.

\medskip The above observation suggests two approaches to fields of definition of the surfaces $X_a$ and their twists : one using
class field theory (see Valloni, \cite{V 21}, \cite{V 23}), the other using the geometry of the moduli space ${\mathcal M}^p_{K3}$  and
elliptic fibrations. This is illustrated by an example due to Brandhorst and Elkies in \cite{BE}; their example
turns out to be isomorphic to a twist of $X_1(7)$ by a prime $O_E$-ideal above $13$  (see Example \ref{BE}).

\medskip
The existence of the $K3$ surfaces is proved by transcendental methods, using the surjectivity of the period map. However, the
method of Brandhorst and Elkies can be used to obtain explicit equations for this family of surfaces. This is done in \cite{BE} for
the above mentioned twist of $X_1(7)$;  I thank Simon Brandhorst
for sending me similar results for two other surfaces in this family, obtained by the method of \cite{BE} (see Examples \ref{7,2} and \ref{2,7}).

\medskip
I thank Simon Brandhorst and Domenico Valloni for interesting discussions, and Matthias Sch\"utt for sending me useful remarks. I am
very grateful to Bill Allombert for his help with PARI/GP.


\section{Lattices, discriminant forms and embeddings}\label{lattice section}

A {\it lattice} is a pair $(L,q)$, where $L$ is a free ${\bf Z}$-module of finite rank, and $q : L \times L \to {\bf Z}$
is a symmetric bilinear form such that ${\rm det}(q) \not = 0$; it is {\it unimodular} if ${\rm det}(1) = \pm 1$, and {\it even} if $q(x,x)$ is an
even integer for all $x \in L$. Set $$L^{\sharp} = \{x \in L \otimes_{\bf Z}{\bf Q} \ | \ q(x,y) \in {\bf Z}  \ {\rm for \ all} \  y \in L \},$$ and
$G_L = L^{\sharp}/L$. The form $q$ induces $G_L \times G_L \to {\bf Q}/{\bf Z}$, called the  {\it discriminant form} of $L$,
and $G_L$ the {\it discriminant group} of $L$; note that the absolute value of ${\rm det}(q)$ is the order of $G_L$.  The
discriminant form is denoted by $(G_L,q_L)$. 

\medskip The Witt group of symmetric bilinear forms on finite abelian groups with values in ${\bf Q}/{\bf Z}$ is denoted by
$W(\bf Q/\bf Z)$; see \cite{Sch}, Chapter V, \S 1. 

\medskip An embedding of lattices $L \to L'$ is called {\it primitive} if its cokernel is free.

\begin{defn} Let $L$ and $L'$ be two lattices. We say that $L$ {\it embeds uniquely into} $L'$ if there exists a primitive embedding
$f : L \to L'$, and if $g : L \to L'$ is another primitive embedding, then there exists $\varphi \in {\rm O}(L')$ such that
$g = \varphi \circ f$. 

\end{defn}


\section{$K3$ surfaces}\label{K3 section basics}

The aim of this section is to recall some basic facts concerning $K3$ surfaces; see \cite{H} or \cite{K} for details.
If $X$ is a complex projective $K3$ surface, set $L_X = H^2(X,{\bf Z})$, and let 
$$H^2(X,{\bf C}) = H^{2,0}(X) \oplus H^{1,1}(X) \oplus H^{0,2}(X)$$ be its Hodge decomposition; we have ${\rm dim}(H^{2,0}) = {\rm dim}(H^{0,2}) = 1$,
and ${\rm dim}(H^{1,1}) = 20$. Let $S_X = L_X \cap H^{1,1}$ be the Picard lattice of $X$, and set $\rho_X = {\rm rank}_{\bf Z}(S_X)$. The
intersection form of $X$ makes $L_X$ into an even unimodular lattice of signature $(3,19)$; since $X$ is projective, the signature of
$S_X$ is $(1,\rho_X -1)$.  Let $T_X$ be the orthogonal complement of $S_X$ in $L_X$. The lattice $T_X$ has signature $(2,20-\rho_X)$, and is
called the {\it transcendental lattice} of $X$. 

\begin{theo}\label{Torelli}  Let $X$ and $Y$ be two complex projective $K3$ surfaces, and let $f : L_X \to L_Y$ be an isometry of lattices whose
$\bf C$-linear extension maps $H^{2,0}(X)$ to $H^{2,0}(Y)$. Then the surfaces $X$ and $Y$ are isomorphic.

\end{theo}

\noindent
{\bf Proof.} This is the weak Torelli theorem, see for instance \cite{H}, Chapter 7, Theorem 5.3. 

\begin{defn} Let $X$ and $Y$ be two complex projective $K3$ surfaces, and let $f : T_X \to T_Y$ be an isometry of lattices. We say
that $f$ is a {\it Hodge isometry} if its $\bf C$-linear extension maps $H^{2,0}(X)$ to $H^{2,0}(Y)$. 

\end{defn}

Let us fix an even unimodular lattice $\Lambda$ of signature $(3,19)$. 

\begin{theo}\label{Hodge theorem}
Let $X$ and $Y$ be two complex projective $K3$ surfaces. Suppose that the lattice $T_X$ embeds uniquely into $\Lambda$, and
that there exists a Hodge isometry $f : T_X \to T_Y$. Then  the surfaces $X$ and $Y$ are isomorphic.
\end{theo}

\noindent
{\bf Proof.} Let us choose isometries $\varphi_X : L_X \to \Lambda$ and $\varphi_Y : L_Y \to \Lambda$; note that $\varphi_X : T_X \to \Lambda$
and $\varphi_Y \circ f : T_X \to \Lambda$ are two primitive embeddings of the lattice $T_X$ into $\Lambda$. Therefore there exists
an isometry $g : L_X \to L_X$ such that $g \circ \varphi_X = \varphi_Y \circ f$; note that the 
$\bf C$-linear extension of
$\varphi_Y^{-1} \circ g \circ \varphi_X : L_X \to L_Y$ sends $H^{2,0}(X)$ to $H^{2,0}(Y)$. Hence by Theorem \ref{Torelli} the $K3$ surfaces
$X$ and $Y$ are isomorphic.

\section{$O_E$-lattices, discriminant ideals  and discriminant modules}\label{discriminant form}

Let $E$ be an algebraic number field with a non-trivial involution $x \mapsto \overline x$, and let $F$ be the fixed field of
this involution; let $n$ be an integer such that $[E:{\bf Q}] = 2n$, and let $\theta \in F^{\times}$ be such that $E = F(\sqrt {\theta})$.  Let $O_E$ be the
ring of integers of $E$, and let $D_{E}$ be the different of $E$.

\medskip
An {\it $O_E$-lattice} is by definition a pair $(I,q)$, where $I$ is a fractional $O_E$-ideal and $q : I \times I \to {\bf Z}$ is given by
$q(x,y) = {\rm Tr}_{E/{\bf Q}}(\alpha x \overline y)$, for some $\alpha \in F^{\times}$; we also use the notation $(I,\alpha)$ 
for this lattice. 

\medskip
If $L$ is an $O_E$-lattice, then so is its dual $L^{\sharp}$, and the quotient $G_L = L^{\sharp}/L$ is also an $O_E$-module, called
the {\it discriminant module} of $L$; it is 
isomorphic to $O_E/{\mathcal D}$ for
some ideal $\mathcal D \subset O_E$. 

\begin{defn} The {\it discriminant ideal} of an $O_E$-lattice $L$ is by definition the ideal $\mathcal D \subset O_E$ such that 
$L^{\sharp}/L$ is isomorphic to $O_E/{\mathcal D}$. The discriminant ideal of $L$ is denoted by $\mathcal D(L)$. 

\end{defn}

\medskip The aim of this section is to characterize the discriminant ideals (equivalently, the discriminant modules) of $O_E$-lattices. 

\medskip
Let $Ram$ be the set of finite places of $F$ that ramify in $E$, and let $Int$ be the set of finite places of $F$ that are inert in $E$. If $v \in Ram$ or $Int$, we denote by $P_v$ the prime $O_E$-ideal corresponding to the
unique place $v_E$ of $E$ above $v$. Let $Ram_E$ be the set of places $w$ of $E$ such that $w = v_E$ for some $v \in Ram$.

\medskip Let $Ram_{odd}$ be the set of  places $v$ of $Ram$ such that $v_E(D_E)$ is odd, and
let $Int_{odd}$ be the set of $v \in Int$ such that $v_E(D_E)$ is odd; let $t$ be the cardinality of $Int_{odd}$.

\medskip
Let us denote by $s$ the number of real embeddings of $F$ that extend to imaginary embeddings of $E$.

\begin{theo}\label{disc} Let $\sigma_1$, $\sigma_2$ be integers $\geqslant 0$ such that $\sigma_1 + \sigma_2 = 2n$. Let $L = (I,q)$ be an $O_E$-lattice
of signature $(\sigma_1,\sigma_2)$.  Then we have

\medskip
{\rm (i)} $\sigma_1 \geqslant n-s$, $\sigma_2 \geqslant n-s$, $\sigma_1  \equiv \ \sigma_2 \equiv \ n-s  \ {\rm (mod \ 2)}$.

\medskip
Let $\mathcal D(L)$ be the discriminant ideal of $L$.  
Then 
$\mathcal D(L) = \prod P^{e_{P}}$ where the product is
taken over the prime ideals of $O_E$ such that the
following conditions hold

\medskip
{\rm (ii)} We have $e_P = 0$ for almost all $P$.

\medskip
{\rm (iii)} For all $P$, we have $e_{\overline P} = e_P$. 

\medskip
{\rm (iv)} If $P = P_v$ with $v \in Ram_{odd}$, then $e_P$ is odd. 

\medskip
{\rm (v)} Let $m$ be the number of $v \in Int$ such that $e_{P_v}$ is odd. If $Ram = \varnothing$, then 
$$\sigma_1 - \sigma_2 \ \equiv 4m \ {\rm (mod \ 8)}.$$

\medskip
Conversely, let $\sigma_1$, $\sigma_2$ be integers $\geqslant 0$ such that $\sigma_1 + \sigma_2 = 2n$ and such that {\rm (i)} holds,
and let $e_P \geqslant 0$ be integers such that $\mathcal D = \prod P^{e_{P}}$ satisfies conditions {\rm (ii)} - {\rm (v)}. Then there exists an $O_E$-lattice $L$ of signature
$(\sigma_1,\sigma_2)$ and discriminant ideal $\mathcal D$. 

\medskip
Moreover, if $a$ is an integer with $0 \leqslant a \leqslant s$ and if $A$ is a set of real places of $F$ of cardinality $a$, then we
can choose $L$ such that $L = (I,\alpha)$ with $\alpha$ negative at the places in $A$ and positive at all the other places of $F$.

\end{theo}

\noindent
{\bf Proof.} (i) follows from \cite{B 99}, Theorem 1, (i). Let $\alpha \in F^{\times}$ be such that $q(x,y) = {\rm Tr}_{E/{\bf Q}}(\alpha x \overline y)$
for all $x,y \in I$. We have $I^{\sharp} = D_{E} ^{-1} \alpha^{-1} \overline I^{-1}$, hence $\mathcal D(L) =  \alpha I \overline I D_E$.
This implies that $\mathcal D(L)$  is of the required form, and that (ii) holds. 
Note that $\overline { D_E} =  D_E$, hence  condition (iii) is satisfied. If  $w = v_E$ for some $v \in Ram_{odd}$,
then $w(\alpha I \overline I D_E)  \equiv \ w( D_E)  \equiv \ 1 \ {\rm (mod \ 2)}$; this implies (iv). 
Finally, let us
prove (v). 
Let $a$ be the number of infinite places $v$ of $F$ such that $(\alpha,\theta)_v = -1$, and let $a'$ be the number of $v \in Int$ such
that $v(\alpha)$ is odd; note that $v(\alpha)$ is odd if and only if $(\alpha,\theta)_v = -1$. Assume that $Ram = \varnothing$; then 
the product formula implies that $a' \equiv a \ {\rm (mod \ 2)}$. We have  $t \equiv a' + m \ {\rm (mod \ 2)}$ by definition, hence
$t \equiv a + m \ {\rm (mod \ 2)}$.

\medskip
We have $\sigma_1 - \sigma_2 = 2s - 4a$ (see for instance \cite {B 99}, Proposition 2.2), and $s \equiv 2t \ {\rm (mod \ 4)}$ (cf. 
\cite{BM}, Theorem 1.6). Hence  $\sigma_1 - \sigma_2 
= 2s - 4a \equiv 4t - 4a \ \equiv \ 4(a + m) - 4a \ \equiv 4m \ {\rm (mod \ 8)}$, as claimed.

\medskip
Conversely, suppose that  conditions (i) - (v) hold. Let $a$ be an integer such that $0 \leqslant a \leqslant s$ and that $\sigma_1 - \sigma_2 =
2s-4a$; such an integer exists by (i). Let $A$ be a set of real places of $F$ of cardinality $a$  that extend to imaginary places of $K$;
this is possible since $a \leqslant s$. 

\medskip Let $M$ be the set of $v \in Int$ such that $e_{P_v}$ is odd; recall that the cardinality of $M$ is denoted by $m$. Let
$A'$ be the symmetric difference of $M$ and $Int_{odd}$, and let $a'$ be the cardinality of $A'$. We define $\epsilon_v = \pm 1$ for
all the places of $v$ as follows. Set $\epsilon_v = -1$ if $v \in A \cup A'$. If $Ram = \varnothing$, then we set $\epsilon_v = 1$
for all the other places of $v$. Otherwise, let us choose a finite place $w$ of $F$ that ramifies in $E$, and set $\epsilon_w =
(-1)^{a + a'}$; set $\epsilon_v = 1$ for all the other places of $F$. 

\medskip
We have $\underset{v} \prod \epsilon_v = 1$.  This is clear if $Ram \not = \varnothing$. Suppose that $Ram = \varnothing$; then 
$\underset{v} \prod \epsilon_v = (-1)^{a+a'}$. Condition (v) implies that $\sigma_1 - \sigma_2 \ \equiv 4m \ {\rm (mod \ 8)}$; since
$\sigma_1 - \sigma_2 = 2s-4a$, this implies that $s - 2a \ \equiv 2m \ {\rm (mod \ 4)}$. We have 
$s \equiv 2t \ {\rm (mod \ 4)}$ by \cite{BM}, Theorem 1.6 and $m \equiv t + a' \ {\rm (mod \ 2)}$ by construction; 
 this implies that $a' \equiv a \ {\rm (mod \ 2)}$, hence 
$\underset{v} \prod \epsilon_v = 1$. 

\medskip There exists $\alpha \in F^{\times}$ such that $(\alpha,\theta)_v = \epsilon_v$ for all places $v$ of $F$ (see for
instance \cite{OM}, Theorem 71.19). We have $v_E(\alpha D_E) \equiv 1 \ {\rm (mod \ 2)}$ if $v \in Ram_{odd}$ or $v \in M$. Let $I$
be an $O_E$-ideal such that $w(\alpha I \overline I D_E) = e_P$ for all places $w$ of $E$, where $P$ is such that $w(P) = 1$. 

\medskip
The $O_E$-lattice $(I,q)$ given by $q(x,y) = {\rm Tr}_{E/{\bf Q}} (\alpha x \overline y)$ has signature $(\sigma_1,\sigma_2)$ 
and discriminant ideal $\mathcal D$. Moreover, $\alpha$ is negative at the places in $A$ and positive at all the other places of $F$. 
This completes the proof of the theorem. 

\begin{coro}\label{disco} Let $\sigma_1$, $\sigma_2$ be integers $\geqslant 0$ such that $\sigma_1 + \sigma_2 = 2n$. Let $L = (I,q)$ be an $O_E$-lattice
of signature $(\sigma_1,\sigma_2)$, and let $G_L$ be the discriminant module of $L$. Then $G_L \simeq \bigoplus O_E /P^{e_{P}}$ where the sum is
taken over the prime ideals of $O_E$ such that conditions {\rm (ii)}-{\rm (v)} above hold.

\medskip
Conversely, let $\sigma_1$, $\sigma_2$ be integers $\geqslant 0$ such that $\sigma_1 + \sigma_2 = 2n$ and such that {\rm (i)} holds,
and let $e_P \geqslant 0$ be integers such that $G = \bigoplus O_E /P^{e_{P}}$ satisfies conditions {\rm (ii)} - {\rm (v)}. Then there exists an $O_E$-lattice $L$ of signature
$(\sigma_1,\sigma_2)$ and discriminant module $G$. 

\medskip
Moreover, if $a$ is an integer with $0 \leqslant a \leqslant s$ and if $A$ is a set of real places of $F$ of cardinality $a$, then we
can choose $L$ such that $L = (I,\alpha)$ with $\alpha$ negative at the places in $A$ and positive at all the other places of $F$.

\end{coro}

\noindent {\bf Proof.} This is an immediate consequence of Theorem \ref{disc}.

\bigskip
The following results will be used in the next sections.

\begin{lemma}\label{even} If no dyadic place of $F$ ramifies in $E$, then every $O_E$-lattice is even. 

\end{lemma}

\noindent
{\bf Proof.} See for instance \cite{B 99}, Proposition 1.

\begin{lemma}\label{even bis} Let $I$ be an ideal of $O_E$, let $\alpha \in F^{\times}$ and let $L = (I,q)$ with 
$q(x,y) = {\rm Tr}_{E/{\bf Q}} (\alpha x \overline y)$. Suppose that for all dyadic places $w$ of $E$ we have
$w(\alpha I \overline I) \geqslant 0$. Then $L$ is an even lattice.

\end{lemma}

\noindent
{\bf Proof.} Let $w$ be a dyadic place of $E$. Since $w(\alpha I \overline I) \geqslant 0$, we have $\alpha x \overline y \in O_w$ for all $x,y \in I$. 
On the other hand, $\overline \alpha = \alpha$, therefore
${\rm Tr}_{K_w/{\bf Q}_2}(\alpha x \overline x)$ is divisible by $2$. This implies that $L$ is even. 

\medskip
Recall that $s$ is the number of real places of $F$ that extend to imaginary places of $E$.

\begin{prop} Suppose that no finite place of $F$ ramifies in $E$. Then $s \equiv 0 \ {\rm (mod \ 2)}$.

\end{prop}

\noindent {\bf Proof.} Let $S$ be the set of real places
of $F$ that extend to imaginary places of $E$; if $v$ is an infinite place of $F$, then $(-1,\theta)_v = -1$ if and only if $v \in S$.
If $v$ is a finite place of $F$, then $v$ is unramified in $E$ by hypothesis, hence $(-1,\theta)_v = 1$. Therefore the product formula implies that $\underset{v \in S} \prod (-1,\theta)_v = 1$, hence $s$ is even, as claimed.

\begin{coro}\label{CM even degree} Assume that $E$ is a CM field with maximal totally real subfield $F$, and that no finite place of $F$ ramifies in $E$. Then $n$
is even.

\end{coro}

\noindent
{\bf Proof.} This follows from the previous proposition, since $s= n$. 

\begin{defn}\label{definition of isomorphic lattices} Let $(L,q)$ and $(L',q')$ be two $O_E$-lattices. We say that $L$ and $L'$ are isomorphic (as $O_E$-lattices)
if there exists an isomorphism of $O_E$-modules $f : L \to L'$ such that $q'(f(x),f(y)) = q(x,y)$ for all $x,y \in L$.

\medskip Set $(V,q) = (L,q) \otimes_{\bf Z} {\bf Q}$ and $(V',q) = (L',q') \otimes_{\bf Z} {\bf Q}$. We say that $L$ and $L'$ {\it become isomorphic
over} $\bf Q$ if there exists an isomorphism of $E$-vector spaces $f : V \to V'$ such that $q'(f(x),f(y)) = q(x,y)$ for all $x,y \in V$.

\end{defn}

\medskip

If $(L,q)$ is an $O_E$-lattice, then there exists an $O_E$-ideal $I$ and $\alpha \in F^{\times}$ such that $q(x,y) = {\rm Tr}_{E/{\bf Q}}(\alpha x \overline y)$. If $L$ is given by $(I,\alpha)$ and $L'$ by $(J,\beta)$ as above, then the $O_E$-lattices $L$ and $L'$ are isomorphic if and only if
there exists $e \in E^{\times}$ such that $J = eI$ and that $\alpha = e \overline e \beta$.

\begin{defn}\label{signature}  Let $L = (I,\alpha)$ and $L' = (J,\beta)$ be two $O_E$-lattices. We say that $L$ and $L'$ have the same signature (as
$O_E$-lattices) if and only if $\tau(\alpha \beta) > 0$ for all embeddings $\tau : F \to {\bf R}$ that extend to imaginary embeddings of $E$. 

\end{defn}

\begin{prop}\label{Q} Let $L = (I,\alpha)$ and $L' = (J,\beta)$ be two $O_E$-lattices of the same signature, and let $G_L = \bigoplus O_E /P^{e_{P}}$,
$G_{L'} =  \bigoplus O_E /P^{e'_{P}}$ 
be their
discriminant modules.  The $O_E$-lattices $L$ and $L'$ become isomorphic over $\bf Q$ if and
only if 
$$e_P  \ \equiv  \ e'_P \ {\rm (mod \ 2)}$$ for all prime $O_E$-ideals $P$ such that $\overline P = P$. 

\end{prop} 

\noindent
{\bf Proof.} Recall that $E = F (\sqrt \theta)$. The $O_E$-lattices $L$ and $L'$ become isomorphic over $\bf Q$ if and only if $\alpha \beta^{-1}$ belongs to ${\rm N}_{E/F} (E^{\times})$;
equivalently, if $(\alpha,\theta)_v = (\beta,\theta)_v$ for all places $v$ of $F$. 
 If $v$ is a real place of $F$, then  $(\alpha,\theta)_v = (\beta,\theta)_v$ if and only if $\alpha \beta$ is positive at $v$,
i.e. if $\tau_v(\alpha \beta) > 0$ where $\tau_v : F \to {\bf R}$ is the embedding corresponding to $v$; since the lattices have
the same signature, the condition holds at $v$. The condition trivially holds if $v$ is an imaginary place of $F$, hence it holds
for all infinite places. Suppose now that $v$ is a finite place of $F$.  If $v$ is split in $E$, then 
$(\alpha,\theta)_v = (\beta,\theta)_v = 1$, hence the condition also holds at split places. If $v$ is inert or ramified in $E$, then $(\alpha,\theta)_v = (\beta,\theta)_v$ if and
only if $e_{P_v}  \ \equiv  \ e'_{P_v} \ {\rm (mod \ 2)}$, where $P_v$ is the unique prime ideal of $O_E$ that is above the prime $O_F$-ideal
corresponding to $v$. This completes the
proof of the proposition.

\begin{notation}
We denote by $h_E$ the class number of $E$. 
\end{notation}

\begin{prop}\label{Simon} Suppose that $E$ is a CM field with maximal totally real subfield $F$,  and that $h_E = 1$. Then two $O_E$-lattices  of the same
signature are isomorphic if and only they if their $O_E$-discriminant modules
are isomorphic.

\end{prop}

\noindent
{\bf Proof.} Let $L$ and $L'$ be two $O_E$-lattices of the same signature. It is clear that if $L$ and $L'$ are isomorphic $O_E$-lattices,
then their $O_E$-discriminant modules are isomorphic.

\medskip
Conversely, suppose that the $O_E$-discriminant modules of $L$ and $L'$ are isomorphic.
Since $h_E = 1$, we have $L \simeq (O_E,\alpha)$ and
$L' \simeq O_E,\beta$ for some 
$\alpha, \beta \in E_0^{\times}$. Since $G_L \simeq G_{L'}$ as $O_E$-modules, the discriminant ideals of $L$ and $L'$ are equal,
and hence $\alpha O_E = \beta O_E$, therefore $\beta \alpha^{-1}$ is a unit of $O_F$. Moreover, it is a totally positive unit, because 
$L$ and $L'$ have the same signature.
Since $h_E = 1$, every totally positive unit
of $E_0$ is a norm of a unit of $E$ (see \cite{Sh}, Proposition A2,  \cite{B 84}, Lemma 3.2). This implies that the $O_E$-lattices
$L$ and $L'$ are isomorphic.

\section{Classification}\label{classification section}

We keep the notation of \S \ref{discriminant form}; in particular, $E$ is a number field of degree $2n$ with a non-trivial involution and
fixed field $F$. To complement the results of \S \ref{discriminant form}, in this section we describe the isomorphism classes of $O_E$-lattices having the same discriminant ideal
and the same signature (see definition \ref{signature}). 

\medskip

Let $C(E)$ be the set of pairs $(I,\alpha)$, where $I$ is an $O_E$-ideal and $\alpha \in F^{\times}$ such that $\alpha I \overline I = O_E$,
and let $C^+(E)$ be the subset of pairs with $\sigma (\alpha) > 0$ for all real embeddings $\sigma$ of $F$ that extend to imaginary embeddings
of $E$. Let us consider the equivalence relation on $C(E)$ (respectively $C^+(E)$) given by

\medskip

\centerline {$(I,\alpha) \equiv (I',\alpha')$ $\iff$ there exists $\gamma \in E^{\times}$ with $I' = \gamma I$ and $\alpha' = ( {1 \over {\gamma \overline \gamma}}) \alpha$.} 

\medskip

In both cases, the
set of equivalence classes is a finite abelian group, the multiplication being induced by $(I,\alpha)(I',\alpha') = (II', \alpha \alpha')$; we denote
this group by $\mathcal C(E)$ (respectively $\mathcal C^+(E)$). 

\medskip Let $L = (J,\beta)$ be an $O_E$-lattice, and let  $a = (I,\alpha) \in C(E)$. Setting $a.L = (IJ, \alpha \beta)$ induces an action of $\mathcal C(E)$ on the set isomorphism 
classes of $O_E$-lattices with the same discriminant ideal as $L$; if moreover, $a \in C^+(E)$, then $a.L$ has the same signature
as $L$, hence we obtain an action of the group $\mathcal C(E)^+$ on the set of isomorphism classes of $O_E$-lattices with the
same discriminant ideal and the same signature as $L$.

\begin{prop}\label{A} 
{\rm (i)} The set of isomorphism classes of $O_E$-lattices with the same discriminant ideal is a principal homogeneous
space over the group $\mathcal C(E)$. 

\medskip
{\rm (ii)} The set of isomorphism classes of $O_E$-lattices with the same discriminant ideal and the same signature is a principal homogeneous
space over the group $\mathcal C^+(E)$. 

\end{prop}

\noindent
{\bf Proof.} Let $L = (J,\beta)$ and $L' = (J',\beta')$ are $O_K$-lattices. We have  
$a.L= L'$ for $a = (J'J^{-1}, \beta' \beta^{-1})$; if $L$ and $L'$ have the same discriminant
ideal, then $a \in C(E)$; if moreover $L$ and $L'$ also have the
same signature, then $a \in C(E)^+$. 

\medskip
Suppose that $E$ is a CM field, and that $F$ is the maximal totally real subfield of $E$. 

\begin{notation} Let ${\mathcal Cl}(E/F)$ 
be the relative class group of $E/F$, and let ${\mathcal Cl}^+(E/F)$ be the strict relative class group. Let $O_F^{\times +}$ be the group
of totally positive units of $O_F$, and let ${\rm N} : E^{\times} \to F^{\times}$ be the norm map. 
\end{notation}

\begin{prop} We have the exact sequences

$$1 \to O^{\times}_F/{\rm N}(O^{\times}_E) \to {\mathcal C}(E) \to {\mathcal Cl}(E/F) \to 1$$

$$1 \to O^{\times +}_F/{\rm N}(O^{\times}_E)  \to {\mathcal C}(E) \to {\mathcal Cl}^+(E/F) \to 1.$$

\end{prop}

\noindent {\bf Proof.} 
The maps ${\mathcal C}(E) \to {\mathcal Cl}(E/F)$ and  ${\mathcal C}(E) \to {\mathcal Cl}^+(E/F)$ are induced
by $(I,\alpha) \mapsto I$. It is easy to check that this gives rise to the above exact sequences.

\section{Twisting}\label{twisting section}

We keep the notation of \S \ref{discriminant form}.

\begin{defn} Let $L$ be an $O_E$-lattice, and let
$J \subset O_E$ be an $O_E$-ideal prime to the discriminant ideal $\mathcal D(L)$ such that $\overline J = J$.
We
say that an $O_E$-lattice $L'$ is a {\it twist of} $L$ {\it by} $J$ if $L$ and $L'$ have the same signature (see Definition \ref{signature}), and 
if  $\mathcal D(L') = \mathcal D(L) J$. 

\end{defn}

We first examine the conditions under which an $O_E$-lattice $L$ has a twist by an ideal $J$.

\begin{prop}\label{existence of twisting} An $O_E$-lattice $L$ has a twist by an ideal $J$ 
if and only if $\mathcal D(L) J$ satisfies condition {\rm (v)} of Theorem \ref{disc}. 

\end{prop}

\noindent
{\bf Proof.} Let $(\sigma_1,\sigma_2)$ the signature of $L$; since $L$ is an $O_E$-lattice, condition
(i) of Theorem \ref{disc} holds. Condition (ii) obviously holds for $\mathcal D(L) J$, and condition (iii) is also satisfied, since $\overline J = J$;
condition (iv) also holds since $J$ is supposed to be prime to $\mathcal D(L')$, hence $J$ does not have any factor $P_v$ with
$v \in Ram_{odd}$. 
Therefore Theorem \ref{disc} implies that there exists an $O_E$-lattice $L'$ with $\mathcal D(L') = \mathcal D(L) J$ if and
only if condition {\rm (v)} hold for $\mathcal D(L) J$. Moreover, Theorem \ref{disc} shows that we can choose
$L'$ with the same signature (as $O_E$-lattice) as $L$. 
 This concludes the proof of the proposition. 
 
 \begin{coro}\label{ramified twist} Suppose that there exists a finite prime of $F$ that ramifies in $E$.
 Let $L$ be an $O_E$-lattice, and let
$J \subset O_E$ be an $O_E$-ideal prime to the discriminant ideal $\mathcal D(L)$ such that $\overline J = J$.
Then $L$ has a twist by $J$.

 \end{coro}
 
 \noindent
 {\bf Proof.} This follows from Proposition \ref{existence of twisting}; indeed, condition (v) is trivially satisfied since $Ram \not = \varnothing$.

\medskip
Such a twist is not unique in general; however, we have the following

\begin{prop}\label{twisting proposition} Suppose that $E$ is a CM field. Let $L$ be an $O_K$-lattice and let $J \subset O_E$ be an ideal
prime to the discriminant ideal $\mathcal D(L)$ such that $\overline J = J$. 
Let $L_1$ and $L_2$ be two twists of $L$ by $J$. If $h_E = 1$, 
then the $O_E$-lattices $L_1$ and $L_2$ are isomorphic.

\end{prop}

\noindent
{\bf Proof.} This follows from Proposition \ref{Simon}. 

\begin{example} Let $L = (I,\alpha)$ be an $O_E$-lattice, and let $P$ be a prime $O_E$-ideal such that $\overline P \not = P$. 
Suppose that $P$ and $\overline P$ are prime to $I$. 
Then $L' = (PI,\alpha)$ is a twist of $L$
by $P \overline P$. Moreover, Proposition \ref{Q} implies that $L$ and $L'$ become isomorphic over $\bf Q$.

\end{example}

\begin{prop} Let $J \subset O_E$ be an ideal with $\overline J = J$.  If all 
prime factors $P$ of $J$ are such that $P \not = \overline P$, then every twist of $O_E$-lattice $L$ by $J$ becomes isomorphic
to $L$ over $\bf Q$.

\end{prop}

\noindent
{\bf Proof.} This is a consequence of Proposition \ref{Q}. 

\begin{example} Let $L $ be an $O_E$-lattice, 
and let $P$ be a prime $O_E$-ideal prime to $\mathcal D(L)$ such that $\overline P  = P$. 
Theorem \ref{disc} and Proposition \ref{existence of twisting} imply that $L$ has a twist by $P$ if and only if $Ram \not = \varnothing$; such a twist does not
necessarily become isomorphic to $L$ over $\bf Q$.

\end{example} 

\section{$O_E$-lattices and $K3$ surfaces}\label{K3 section}

We keep the notation of  the previous sections, and set $[E:{\bf Q}] = 2n$. We assume in addition that $E$ is a CM field with $[E:{\bf Q}] \leqslant 20$, i.e. $E$ is totally imaginary and
$F$ is totally real with  $1 \leqslant n \leqslant 10$.  If $X$ is a $K3$ surface,
we denote by $T_X$ its transcendental lattice.

\medskip
Let us fix an even unimodular lattice $\Lambda$ of signature $(3,19)$.

\begin{prop}\label{existence}  Let $L$ be an $O_E$-lattice of signature $(2,2n-2)$ and assume that $L$ embeds primitively into $\Lambda$. Then there exists
a complex projective $K3$ surface $X$ such that the transcendental lattice $T_X$ of $X$ has a structure of $O_E$-lattice isomorphic to $L$,
and ${\rm End}_{\rm Hdg}(T_X) = O_E$. 

\end{prop} 

\noindent
{\bf Proof.} Let $L$ be given by $(I,\alpha)$, where $I$ is an $O_E$-ideal and $\alpha \in F^{\times}$. Let $\sigma' : F \to {\bf R}$ be the
real embedding of $F$ such that $\sigma'(\alpha) > 0$; note that $\alpha$ is negative at all the other real embeddings of $F$, since the signature of $L$ is $(2,2n-2)$. Let $\sigma : E \to {\bf C}$ be an extension of $\sigma'$ to $E$. 

\medskip
Let $f : L \to \Lambda$ be a primitive embedding, and let us also denote by $f$ its extension to $f : L \otimes_{\bf Z}{\bf C} \to \Lambda \otimes_{\bf Z}{\bf C}$. We have $ L \otimes_{\bf Z}{\bf C} = \underset{\tau : E \to {\bf C}} \oplus {\bf C}_{\tau}$. Set $\Lambda^{2,0} = {\bf C}_{\sigma}$. 

\medskip
We obtain the desired $K3$ surface by surjectivity of the period map. Indeed, the choice of $\Lambda^{2,0}$ induces on $\Lambda$ a Hodge
structure. Let $X$ be the corresponding $K3$ surface. By construction, we have $T_X = L$ and ${\rm End}_{\rm Hdg}(T_X) = O_E$.
Since the signature of $L$ is $(2,2n-2)$, the surface $X$ is projective.

\begin{notation} If $X$ is a $K3$ surface with complex multiplication by $O_E$, we denote by $G_X$ the discriminant $O_E$-module $T_X^{\sharp}/T_X$. 
The minimal number of generators (as an abelian group) of $G_X$ is denoted by $\ell(X)$, and is called the {\it length} of $X$. 
\end{notation}

\begin{defn}\label{K3 twist} Let $X$ and $Y$ be two complex projective $K3$ surfaces with complex multiplication by $O_E$. Let $J \subset O_E$ be
an $O_E$-ideal such that $\overline J = J$. We say that $Y$ is a {\it twist of $X$ by $J$} if $G_{Y} \simeq G_X \oplus O_E/J$. 

\end{defn}

Let us consider $E$ embedded in $\bf C$. This implies that if $X$ and $Y$ are two $K3$ surfaces with CM by $O_E$, then the
$O_E$-lattices $T_X$ and $T_Y$ have the same signature (as $O_E$-lattices).

\begin{prop}\label{iso surfaces} Let $L$ be an $O_E$-lattice, and assume that $L$ embeds uniquely into $\Lambda$. 
Let $X$ and $Y$ be two complex $K3$ surfaces with CM by $O_E$, and suppose that the $O_E$-lattices $T_X$ and $T_Y$ are
isomorphic to $L$. Then the surfaces $X$ and $Y$ are isomorphic.

\end{prop} 

\noindent
{\bf Proof.} This follows from Theorem \ref{Hodge theorem}.

\medskip

We next note that $h_E = 1$, then $K3$ surfaces with maximal complex multiplication by $E$ of length $\leqslant 20-2n$
 are determined by their discriminant modules.

\begin{prop}\label{Nikulin} Suppose that $h_E = 1$, and let $X$ and $Y$ be two $K3$ surfaces with maximal complex multiplication by $E$ of length 
$\leqslant 20-2n$.
Then $X$ and $Y$ are isomorphic
if and only if the discriminant $O_E$-modules of $T_X$ and $T_{Y}$ are isomorphic.

\end{prop}

\noindent
{\bf Proof.} if $X$ and $Y$ are isomorphic, then the $O_E$-lattices $T_X$, $T_Y$ are isomorphic, and hence so are their
discriminant modules.  Let us prove the converse. Since $h_E = 1$, 
Proposition \ref{Simon} implies that the 
$O_E$-lattices $T_X$ and $T_Y$ are isomorphic. 

\medskip
We fix an even unimodular lattice $\Lambda$ of signature $(3,19)$.
Since by hypothesis $\ell(X), \ell(Y)  \leqslant 20 - 2n$, the lattices $T_X$, $T_Y$ are uniquely embedded in $\Lambda$ (see Nikulin \cite{N}, 
Theorem 1.14.4).  Therefore by Proposition \ref{iso surfaces} the $K3$ surfaces X and $Y$ are isomorphic. 

\begin{coro}\label{twist uniqueness} Let $X$ be a complex $K3$ surface with CM by $O_E$, and suppose that $h_E = 1$. Let 
Let $J \subset O_E$ be
an $O_E$-ideal such that $\overline J = J$, and let $Y_1$, $Y_2$ be two twists of $X$ by $J$ of length $\leqslant 20 - 2n$. Then $Y_1$ and $Y_2$ are
isomorphic. 

\end{coro}

\noindent
{\bf Proof.} This follows from Proposition \ref{Nikulin}.

\section{Existence of $K3$ surfaces with maximal complex multiplication}\label{CM section}

We keep the notation of the previous sections; in particular, $E$ is a CM field of degree $\leqslant 20$. The aim of this
section is to give a criterion for the existence of infinitely many isomorphism classes of $K3$ surfaces with complex multiplication by $O_E$. Valloni proved that
this is always the case if the degree of $E$ is $\leqslant 10$ (cf. \cite{V 21}, Proposition 6.11); as we will see, this result extends to fields of degree $\leqslant 14$
(see Corollary \ref{14}).

\medskip
We start
by introducing some notation.
If $w$ is a place of $E$, we denote by $f_w$ its residual degree.

\begin{notation} If $p$ is a prime number such that $p \not = 2$, we denote by $Ram(p)$ the set of places of $Ram_E$ above $p$. Let
$Ram(2)$ be the set of dyadic places of $E$ such that $w(D_E) > 0$. For all prime numbers $p$, set

$$f(p) = \underset{w \in Ram(p)} \sum f_w.$$

\end{notation}

Note that for almost all $p$, we have $Ram(p) = \varnothing$, hence $f(p) = 0$. 

\begin{theo}\label{CM theorem} Suppose that 
$f(p) < 22 - 2n$ for all prime numbers $p$ such that $Ram(p) \not = \varnothing$. 
Then there exist infinitely many non-isomorphic complex projective $K3$ surfaces having complex multiplication by $O_E$. 

\medskip
Moreover, there exist infinitely many such surfaces in the same isogeny class.

\end{theo}

\noindent {\bf Proof.} Suppose first that $Ram \not = \varnothing$. Let $P$ be a prime  ideal of $O_E$ of degree 1; there exist
infinitely many such ideals by Chebotarev's density theorem. Assume that $P$ is not dyadic, and that ${\rm N}(P)$ is relatively prime to 
${\rm N}(P_v)$ for all $v \in Ram$, where ${\rm N}$ is the norm map. 

\medskip
By
Corollary  \ref{disco}, there exists an $O_E$-lattice with signature $(2,2n-2)$ and discriminant module $\underset{v \in Ram_{odd}} \oplus O_E/P_v \oplus O_E/P \oplus O_E/{\overline P}$; let $I$ be an $O_E$-ideal
and let $\alpha \in F^{\times}$ such that the lattice $q : I \times I \to {\bf Z}$ with  $q(x,y) = {\rm Tr}_{E/{\bf Q}}(\alpha x \overline y)$ for all $x,y \in I$
is such a lattice. 

\medskip
If there exist dyadic places of $F$ that ramify in $E$, 
let $J = I \underset{w \in Ram(2)} \prod P_w^{e_w}$ with $e_w \in {\bf Z}$ 
such that $w(\alpha J \overline J) \geqslant 0$ for all dyadic places $w$ of $E$, where $P_w$ is the prime $O_E$-ideal such that $w(P_w) = 1$. 
Let $L$ be the $O_E$-lattice given by $q : J \times J \to {\bf Z}$ such that $q(x,y) = {\rm Tr}_{E/{\bf Q}}(\alpha x \overline y)$ for all $x,y \in J$.
By Lemma \ref{even} and Lemma \ref{even bis}, the lattice $L$ is even. 

\medskip 
Assume now that $Ram = \varnothing$, i.e. no finite place of $F$ ramifies in $E$;  by Lemma \ref{CM even degree} this implies that $n$ is even.  As before, let $P$ be a non-dyadic prime ideal of $O_E$ of degree $1$. 
If $n \equiv 2\ {\rm (mod \ 4)}$, set $G = O_E/P  \oplus O_E/{\overline P}$. Suppose that $n \equiv 0 \  {\rm (mod \ 4)}$. Since $E$ is a CM field, Chebotarev's density
theorem implies that there exist infinitely many prime $O_E$-ideals $Q$ such that $Q \cap O_F$ is inert in $E/F$, and that the residual degree
of $Q$ is 2. Let $Q$ be such an ideal, and set $G = O_E/P  \oplus O_E/{\overline P}  \oplus O_E/Q$ in this case. By Corollary \ref {disco}, there exists an $O_E$-lattice $L$
of signature $(2,2n-2)$ and discriminant module $G$; by Lemma \ref{even}, this lattice is even.

\medskip In all the above cases, we denote by $G$ the discriminant module of $L$, and let $\ell(G)$ be the number of generators of $G$
as an abelian group. 
Let 
$f$ be the maximum of the integers $f(p)$ such that $Ram(p) \not = \varnothing$ if there exists such a $p$ with $f(p) > 1$; otherwise,
set $f = 2$. We have $\ell(G) \leqslant f$. 
If $n \leqslant 9$, then the hypothesis
implies that $\ell(G) < 22 - 2n$; by Nikulin's result \cite{N}, Corollary 1.12.3 this implies that
$L$ can be primitively embedded in an even, unimodular lattice $\Lambda$ of signature $(3,19)$. Suppose that $n = 10$, and
let $p = N(P)$. The $p$-component of $G$ is ${(\bf F}_p)^2$, and the $p$-component of the discriminant form has determinant $-p^2$.
By \cite{N}, Theorem 1.12.2, the lattice $L$ can be primitively embedded in an even, unimodular lattice $\Lambda$ of signature $(3,19)$ 
in this case as well. 

\medskip 
Let $\sigma' : F \to {\bf R}$ be the unique embedding of $F$ such that $\sigma'(\alpha) > 0$, and let $\sigma : E \to {\bf C}$ be one
of the two extensions of $\sigma'$ to $E$. Let $J \otimes_{\bf Z} {\bf C} = \underset {\tau : E \to {\bf C}} \oplus {\bf C}_{\tau}$ 
and set $\Lambda^{2,0} = {\bf C}_{\sigma}$, where we consider $J \otimes_{\bf Z} {\bf C}$ contained in $\Lambda \otimes_{\bf Z} {\bf C}$.
This endows the lattice $\Lambda$ with a Hodge structure. Let $X$ be the corresponding $K3$ surface : such a surface exists by the
surjectivity of the period map. The transcendental lattice of $X$ is isomorphic to $L$, a lattice of signature $(2,2n-2)$, hence the
surface $X$ is projective. It has complex multiplication by $O_E$ by construction. Varying the ideal $P$ gives rise to infinitely
many non-isomorphic projective $K3$ surfaces having complex multiplication by $O_E$. 
By Proposition \ref{Q} the $O_E$-lattices become isomorphic over $\bf Q$, hence the $K3$ surfaces are all isogeneous (cf. \cite{Mu}, \cite {N 87},
\cite {Bu}). This completes the proof
of the theorem. 

\bigskip

We now state some consequences of this result.

\begin{coro}\label{unramified} If no finite place of $F$ ramifies in $E$,  then there exist infinitely many non-isomorphic complex projective  $K3$ surfaces with
complex multiplication by $O_E$.

\end{coro}

\noindent
{\bf Proof.} $Ram = \varnothing$ in this case, hence $f(p) = 0$ for all prime numbers $p$.

\begin{coro}\label{14} If $[E:{\bf Q}] \leqslant 14$, then there exist infinitely many non-isomorphic complex projective  $K3$ surfaces with
complex multiplication by $O_E$.

\end{coro}

\noindent
{\bf Proof.} Recall that $[E:{\bf Q}] = 2n$.  It is easy to see that for all prime numbers $p$, we have $f(p) < n$. 
If $[E:{\bf Q}] \leqslant 14$, then $f(p) < 7$, and hence $f(p) < 22 - 2n \leqslant 8$. Therefore by Theorem \ref{CM theorem} there 
exist infinitely many non-isomorphic complex projective $K3$ surfaces with complex multiplication by $O_E$. 

\medskip Note that this no longer holds in general when $[E:{\bf Q}] \geqslant 16$, as shown by the following proposition and example.

\begin{prop}\label{high ramification} Let $p$ be a prime number, $p \not = 2$. Suppose that there exists a prime $O_F$-ideal $P$
that ramifies in $E$ such that $N(P) = p^f$ with $f > 22-2n$. Then there does not exist any complex projective  $K3$ surfaces with
complex multiplication by $O_E$.

\end{prop}

\noindent
{\bf Proof.} Let $X$ be a complex projective $K3$ surface,  let $T_X$ be the transcendental lattice and let $S_X$ be the Picard lattice of $X$.
If $X$ has complex multiplication by $E$, then ${\rm rank}(T_X) = 2n$ and hence ${\rm rank}(S_X) = 22-2n$. Let $G$ be the discriminant
group of the lattice $T_X$; then $G$ is also the discriminant group of the lattice $S_X$. Theorem \ref{disc} implies that $O_E/P$ is
a subgroup of $G$; the hypothesis on $P$ implies that the minimal number of generators of $G$ is $> 22-2n$. This contradicts the
fact that $G$ is the discriminant group of the lattice $S_X$, of rank $22-2n$. 

\begin{example}\label{counterexample} Let $F$ be the maximal totally real subfield of the cyclotomic field ${\bf Q}(\zeta_{17})$, and
set $E = F(\sqrt {-3})$. Note that $E$ is a CM field of degree $16$. There exists a unique prime ideal $P$ above $3$ in $F$; this
ideal ramifies in $E$, and its residual degree is $8$, i.e. $N(P) = 3^8$; by Proposition \ref{high ramification} this implies that there
does not exist any complex projective  $K3$ surfaces with
complex multiplication by $O_E$.
The same method gives rise to  infinitely many examples in degrees $16$, $18$ and $20$. 

\end{example}

\begin{coro}\label{cyclotomic}  If $E$ is a cyclotomic field with $2 \leqslant [E:{\bf Q}] \leqslant 20$, then there exist infinitely many non-isomorphic complex projective  $K3$ surfaces with
complex multiplication by $O_E$.

\end{coro}

\noindent
{\bf Proof.} If no finite prime of $F$ ramifies in $E$, then this follows from Corollary \ref {unramified}. Suppose now that there exist finite primes of $F$ that
ramifiy in $E$; this implies that $E = {\bf Q}(\zeta_{p^r})$, where
$p$ is a prime number and $r \geqslant 1$ is an integer. Let $P$ be the unique ramified ideal of $O_E$. Then the residual degree of $P$ is $1$,
and $P$
is the only prime ideal of $O_E$ above $p$, hence $f(p) = 1$. We have $f(q) = 0$ for all prime numbers $q \not = p$, hence 
by Theorem \ref{CM theorem}, this implies that there exist infinitely many non-isomorphic complex projective  $K3$ surfaces with
complex multiplication by $O_E$.

\section{$K3$ surfaces with a given discriminant ideal}\label{given section}

We keep the notation of the previous sections : $E$ is a $CM$ field of degree $2n$, with maximal totally real subfield $F$, and $2n \leqslant 20$. 
We now apply the results of \S \ref{discriminant form} - \S \ref{K3  section} to the existence and classification of $K3$ surfaces with
a given discriminant ideal.

\begin{defn} Let $X$ be a $K3$ surface with CM by $O_E$. The {\it discriminant ideal} of $X$, denoted by $\mathcal D_X$, is the
integral ideal of $O_E$ such that the $O_E$-modules $T_X^{\sharp}/T_X$ and $O_E/{\mathcal D}_X$ are isomorphic.

\end{defn}

\medskip

Recall that $G_X = T_X^{\sharp}/T_X$ is called the discriminant module of $X$, and  that the length of $X$, denoted by $\ell(X)$, is
by definition the minimal number of generators of $G_X$, as an abelian group.

\medskip
If $X$ is a $K3$ surface with CM by $O_E$, then $T_X$ is an even $O_E$-lattice of signature $(2,2n-2)$, hence the disciminant ideal  $\mathcal D_X$
satisfies the conditions of Theorem \ref{disc} for $\sigma_1 = 2$ and $\sigma_2 = 2n-2$; moreover, $\ell(X) \leqslant 22-2n$.

\begin{coro}\label{D existence} Suppose that no dyadic place of $F$ ramifies in $E$, and let $\mathcal D \subset O_E$ be an ideal satisfying conditions {\rm (ii) - (v)} of Theorem \ref{disc} for $(\sigma_1,\sigma_2) = (2,2n-2)$,
and suppose that
$\ell(O_E/{\mathcal D}) < 22-2n$. Then there exists a $K3$ surface $X$ with $\mathcal D_X = \mathcal D$. 
\end{coro}

\noindent {\bf Proof.} Theorem \ref{disc} implies that there exists an $O_E$-lattice $L$ with discriminant ideal $\mathcal D$ and signature
$(2,2n-2)$. Since no dyadic place of $F$ ramifies in $E$, the lattice is even (see Lemma \ref{even}). 
The lattice $L$ embeds primitively into
the $K3$-lattice $\Lambda$ (see \cite{N}, Corollary 1.12.3),  hence by Proposition \ref{existence} there exists a complex projective $K3$ surface $X$ with complex
multiplication by $O_E$ such that $T_X \simeq L$. 

\medskip Let us consider $E$ embedded in $\bf C$; hence all transcendental lattices of $K3$ surfaces with maximal
complex multiplication by $E$ have the same signature (as $O_E$-lattices). 

\medskip Set $\mathcal C = \mathcal C^+(E)$, with the notation of \S \ref{classification section}. Let $\mathcal D \subset O_E$ be
an ideal such that $\ell(O_E/{\mathcal D}) < 22-2n$. 

\begin{coro}\label{K3 with D} Suppose that no dyadic place of $F$ ramifies in $E$.  Then the set of complex projective $K3$ surfaces with complex multiplication by $O_E$ and discriminant
ideal $\mathcal D$ is a principal homogeneous space over $\mathcal C$. 

\end{coro} 

\noindent
{\bf Proof.} The hypotheses imply there exists a $K3$ surface with $\mathcal D_X = \mathcal D$
(see Corollary \ref{D existence}), and that the transcendental lattice $T_X$ embeds uniquely into the $K3$-lattice. The
corollary now follows from Proposition \ref{A} (ii), the fact that by Lemma \ref{even} every $O_E$-lattice is even, combined
with Proposition \ref {existence} and Proposition \ref{iso surfaces}. 

\begin{coro}\label{finite} There exist only finitely many isomorphism classes of $K3$ surfaces with $CM$ by $O_E$ and discriminant
ideal $\mathcal D$. 

\end{coro}

\noindent
{\bf Proof.} Recall that $X$ is a $K3$ surface with $CM$ by $O_E$ and discriminant
ideal $\mathcal D$. 
Let $Y$ be another $K3$ surface with these properties. Both $T_X$ and $T_Y$ are even $O_E$-lattices of the same signature, and their discriminant ideals
are equal by hypothesis. Proposition \ref{A} (ii)  implies that there exists $a \in \mathcal C$ such that $a.T_X = T_Y$. The group $\mathcal C$
is finite, hence there are only finitely many possibilities for the isomorphism class of the $O_E$-lattice $T_Y$. By Proposition 
\ref{iso surfaces}, this implies that there are only finitely many isomorphism classes of $K3$ surfaces $Y$ as above. 

\begin{remark} Let $|\mathcal C|$ be the order of the group $\mathcal C$. The number of isomorphism classes of 
$K3$ surfaces with $CM$ by $O_E$ and discriminant
ideal $\mathcal D$ is $\leqslant |\mathcal C|$, and equality holds if no dyadic prime of $F$ ramifies in $E$.

\end{remark}

\section{Picard lattices and complex multiplication}\label{Picard section}
 
 The aim of this section is to discuss the relationship between complex multiplication by a ring of integers, and properties of
 the Picard and transcendental lattices. We keep the notation of the previous sections; in particular, $E$ is a CM field, $F$ is its maximal totally
 real subfield, ${\rm deg}(E) = 2n$, with $2n \leqslant 20$. In this section, we assume that no dyadic prime of $F$ ramifies in $E$.
 
 \begin{prop}\label{unramified} Let $X$ be a complex projective $K3$ surface with maximal complex multiplication by $E$. Suppose
 that the Picard lattice $S_X$ is unimodular {\rm (equivalently, $T_X$ is unimodular)}. Then no finite prime of $F$ ramifies in $E$, and
 $2n \ \equiv 4 \ {\rm (mod \ 8)}.$
 
 \medskip
 Conversely, if no finite prime of $F$ ramifies in $E$ and
 $2n \ \equiv 4 \ {\rm (mod \ 8)}$, then there exists a complex projective $K3$ surface $X$ such that $S_X$ and $T_X$ are unimodular.

 \end{prop}

\noindent
{\bf Proof.} The hypothesis implies that $T_X$ is a unimodular $O_E$-lattice; recall that this lattice is even. With the notation of \S \ref{discriminant form}, 
this implies that $e_{P_v} = 0$ for all places $v$ of $E$, therefore by Theorem \ref{disc} (iv) we have $Ram = \varnothing$. 
We have $(\sigma_1,\sigma_2) = (2,2n-2)$; since $T_X$ is unimodular, $m = 0$. Therefore Theorem \ref{disc} (v) implies that $2n - 4\ \equiv 0 \ {\rm (mod \ 8)}$, hence  $2n \ \equiv 4 \ {\rm (mod \ 8)}$, as claimed. 

\medskip
Conversely, assume that no finite prime of $F$ ramifies in $E$ and that
 $2n \ \equiv 4 \ {\rm (mod \ 8)}$. Then by Theorem \ref{disc} there exists a unimodular $O_E$-lattice $L$ of signature $(2,2n-2)$. 
 Since no finite prime of $F$ ramifies in $E$, this lattice is even (cf. Lemma \ref{even}). By
 Proposition \ref{existence} there exists a complex projective $K3$ surface $X$ with CM by $O_E$ and $T_X \simeq L$. Since
 $S_X$ is the orthogonal complement of $T_X$ in $H^2(X,{\bf Z})$, the lattice $S_X$ is also unimodular. 
 
 \bigskip
 Suppose now that $2n = 20$; in this case, the Picard lattice is of rank 2. We denote by $U$ the rank $2$ hyperbolic lattice, and
 if $N$ is an integer, we denote by $U(N)$ the lattice $U$ with values multiplied by $N$. 
 
 \begin{notation} Let $S_1$ be the set of prime numbers $p$ such that there exists a prime $O_E$-ideal $P$ 
 with $\overline P \not = P$ and  $p =  {\rm N}_{E/{\bf Q}}(P)$, let $S_2$ be the set of prime numbers $p$
such that there exists a prime $O_E$-ideal $P$ with $\overline P = P$ and $p^2 = {\rm N}_{E/{\bf Q}}(P)$, and
let $S_3$ be the set of prime numbers $p$ such that there exists a prime ideal $P$ of $E$ such that $\overline P = P$
and $p =  {\rm N}_{E/{\bf Q}}(P)$. 

 \end{notation}
 
 \begin{lemma} {\rm (i)} The set $S_3$ is finite. {\rm (ii)} If no finite prime of $F$ ramifies in $E$, then $S_3 = \varnothing$. 
 
 \end{lemma}
 
 \noindent
 {\bf Proof.} If $p \in S_3$, then there exists a prime $O_F$-ideal above $p$ that ramifies in $E$; this proves both (i) and (ii).

 \begin{notation}
 Let $\mathcal N_E$ be the set of integers $N \geqslant 1$ such that 
 $N =  \underset {i \in I} \prod p_i^{n_i}$, where for all $i \in I$ we have $p_i \in S_1$ or $p_i \in S_2$,  and $n_i \geqslant 0$
  is an integer such that  if no finite prime of $F$ ramifies in $E$, then $\underset{p_i \in S_2} \sum n_i$ is even. 
 
 \end{notation}

 Recall that if $P$ is a prime ideal of $O_E$ above the prime number $p$, we denote by $f_P$ the residual degree of $P$, i.e. $f_P = [O_E/P : {\bf F}_p]$.
 
 \begin{prop}\label{Picard unramified} Suppose that $[E:{\bf Q}] = 20$, and that $S_3 = \varnothing$. 
  Let $X$ be a complex projective $K3$ surface with maximal complex multiplication by $E$.  
  
  \medskip
  Then the Picard lattice $S_X$ is
  isomorphic to $U(N)$ with $N \in \mathcal N_E$. 
  
\medskip
Conversely, if $N \in \mathcal N_E$,  then there exists a complex projective $K3$ surface with maximal complex multiplication by $E$
with Picard lattice isomorphic to $U(N)$. 
 
 \end{prop} 
 
 \noindent
 {\bf Proof.} The discriminant module $G_X$  is isomorphic to $\underset{P} \bigoplus \ O_E/P^{e_P}$ for some
 prime $O_E$ ideals $P$ and integers $e_P \geqslant 0$ with $e_P = e_{\overline P}$. Since $S_3 = \varnothing$,  the order
 of $G_X$ is a square, and therefore ${\rm det}(S_X)$ is a square. This implies that $S_X \simeq U(N)$ for some integer $N$, and
 $G_X \simeq ({\bf Z}/N {\bf Z})^2$; therefore $\ell(X) = 2$. 
 
 \medskip
 If $P$ is a prime $O_E$-ideal such that $e_P \not = 0$ and $\overline P \not = P$, then  this implies that $f_P = 1$. Set 
  $p = {\rm N}_{E/{\bf Q}}(P)$; we have $p \in S_1$.
  
  \medskip If $P$ is a prime $O_E$-ideal such that $e_P \not = 0$ and $\overline P = P$, then we have $p = {\rm N}_{E/{\bf Q}}(P)$
  or  $p^2 = {\rm N}_{E/{\bf Q}}(P)$, and this implies $p \in S_3$ in the first case and $p \in S_2$ in the second one. But
  $S_3 = \varnothing$ by hypothesis, hence we have $p^2 = {\rm N}_{E/{\bf Q}}(P)$ and $p \in S_2$.
  
  \medskip Note that if a prime number $p$ divides $N$, then there exists a prime $O_E$-ideal above $p$ with $e_P \not = 0$, hence
  we proved that $N$ is a product of primes in $S_1 \cup S_2$. 
  
  \medskip
 Set $N =  \underset {i \in I} \prod p_i^{n_i}$; 
 it remains to prove that if no finite prime of $F$ ramifies in $E$, then $\underset{p_i \in S_2} \sum n_i$ is even.

 \medskip Suppose that no finite prime of $F$ ramifies in $E$, and
 note that under this hypothesis, if $P$ is a prime ideal with $\overline P = P$, then $e_P > 0$ $\iff$ ${\rm N}(P) = p_i^2$ for some $p_i \in S_2$. Moreover,
 $e_P = n_i$. With the notation of Theorem \ref{disc}, we have $\sigma_1 = 2$ and $\sigma_2 = 18$, hence
 $\sigma_1 - \sigma_2 = -16$; by Theorem \ref{disc} (v) this implies that $m$ is even. Therefore
 $e_P$ is odd for an even number of prime $O_E$-ideals with $\overline P = P$;  this implies that the sum
 $\underset{p_i \in S_2} \sum n_i$ is even, as claimed.

 \medskip
 Conversely, let $p = p_i$ be a divisor of $N$, and let $P$ be a prime $O_E$-ideal with $p = {\rm N}_{E/{\bf Q}}(P)$ if $p \in S_1$,
 and $p^2 = {\rm N}_{E/{\bf Q}}(P)$ if $p \in S_2$. Set $e_P = n_i$ and $G = \underset{P} \bigoplus \ O_E/P^{e_P}$. If no finite
 prime of $F$ ramifies in $E$, then 
 $\underset{p_i \in S_2} \sum n_i$ is even, hence the number of prime ideals $P$ with $\overline P = P$ and $e_P > 0$ is even; with
 the notation of Theorem \ref{disc}, this implies that $m$ is even. 
 By Theorem \ref{disc} there
 exists an $O_E$-lattice $T$ of signature $(2,18)$ and discriminant module $G$. The lattice $T$ embeds primitively into
 the $K3$-lattice $\Lambda$, hence by Proposition \ref{existence} there exists a complex projective $K3$ surface $X$ with transcendental
 lattice $T$ and maximal complex multiplication by $O_E$. The orthogonal complement of $T$ in $\Lambda$ is isomorphic
 to $U(N)$ by construction, and this completes the proof of the proposition. 
 
 \begin{remark} If $E$ is a Galois extension of ${\bf Q}$ and if no finite prime of $F$ ramifies in $E$, then $S_1$ is the set of prime numbers that split completely in $E$,
 and $S_2$ is the set of those that split completely in $F$, but not in $E$.
 
 \end{remark}
 
\begin{example}\label {44 and 66}  Let $E = {\bf Q}(\zeta_m)$ with $m = 44$ or $66$, and let $p$ be a prime number. No finite prime of $F$
ramifies in $E$, hence $S_3 = \varnothing$. We have

\medskip

$\bullet$ $p \in S_1 \iff p \equiv 1 \ {\rm (mod \ m)}$; 

\medskip
$\bullet$ $p \in S_2 \iff p^2 \equiv 1 \ {\rm (mod \ m)}$.

\medskip Proposition \ref{Picard unramified}  implies that there exists a $K3$ surface with maximal complex multiplication by $E$ 
with Picard lattice $L$ $\iff$ $L \simeq U(N)$ where $N \geqslant 1$ is an integer  $\equiv 1 \ {\rm (mod \ m)}$
such that  all the prime divisors of $N$ are  $ \equiv \pm 1 \ {\rm (mod \ m)}$.

\medskip Note that for $N = 1$ we recover Kondo's $K3$ surfaces, see \cite{Ko} and \cite{LSY}.

\end{example}

\medskip

Still supposing that $2n = 20$, we now deal with the case where $S_3 \not = \varnothing$. We start by introducing some notation.

\begin{notation} Let $K$ be a real quadratic field, and let
$\sigma : K \to K$ be the unique non-trivial element of the Galois group of $K$ over ${\bf Q}$. Let $O$ be
an order of $K$ and let $I \subset O$ be a projective ideal of $O$. Let ${\rm N} : K \to {\bf Q}$ be the norm map. We denote
by $q_I$ the quadratic form $q_I : I \times I \to {\bf Z}$ defined by $q_I(x,y) = {1 \over {{\rm N}(I)}} {\rm Tr}_{K/{\bf Q}} x \sigma(y)$.

\medskip

The {\it conductor} of an order $O$ is by definition the index
of $O$ in the ring of integers of $K$; we denote by $cond(O)$ the conductor of $O$. 

\end{notation}

\begin{notation} 
 Let $\mathcal M_E$ be the set of integers $N \geqslant 1$ such that 
 $N =  \underset {i \in I} \prod p_i^{n_i}$, where for all $i \in I$ we have $p_i \in S_1$, $p_i \in S_2$ or $p_i \in S_3$,  and $n_i \geqslant 0$
  is an integer such that $n_i$ is even if $p_i \in S_3$. 

\end{notation}

\begin{prop}\label{Picard ramified} Suppose that $[E:{\bf Q}] = 20$ and that $S_3 \not = \varnothing$. 
Let $X$ be a complex projective $K3$ surface with maximal complex multiplication by $E$.

\medskip
{\rm (i)} If ${\rm det}(T_X)$ is a square, then $S_X \simeq U(N)$ for some $N \in \mathcal M_E$.

\medskip
{\rm (ii)} Suppose that ${\rm det}(T_X)$ is not a square, and let ${\rm det}(T_X) = d c^2$, where $d$ is a square-free integer. Set
$K = {\bf Q}(\sqrt d)$. Then 

$$S_X \simeq q_I$$ where $I$
is an $O$-ideal of an order $O$ of $K$ of conductor $c$; moreover, $c \in \mathcal M_E$.

\end{prop}

\noindent
{\bf Proof.} (i) Since ${\rm det}(T_X)$ is a square, $|{\rm det}(S_X)|$ is also a square, hence $S_X \simeq U(N)$ for some
integer $N$. We have $\ell(X) = 2$. The discriminant module $G_X$  is isomorphic to $\underset{P} \bigoplus \ O_E/P^{e_P}$ for some
 prime $O_E$ ideals $P$ and integers $e_P \geqslant 0$ with $e_P = e_{\overline P}$. If a prime number $p$ divides $N$, then
 there exists a prime $O_E$-ideal $P$ above $p$ such that $e_P \not = 0$.
 
 \medskip
 Suppose that $p$ is a prime divisor of $N$ and that $P$ is a prime $O_E$-ideal above $p$ with $e_P \not = 0$ such that $\overline P \not = P$. 
 Since $\ell(X) = 2$, this implies that $f_P = 1$ and $p = {\rm N}_{E/{\bf Q}}(P)$, hence $p \in S_1$.
 
 \medskip
 Let $p$ be a prime divisor of $N$ such that there exists a prime $O_E$-ideal above $p$ with $e_P \not = 0$ and $\overline P = P$. 
 If $f_P = 2$, then $p^2 = {\rm N}_{E/{\bf Q}}(P)$, and $p \in S_2$. Suppose that $f_P = 1$. Then we have $p = {\rm N}_{E/{\bf Q}}(P)$,
 hence $p \in S_3$. Since ${\rm det}(T_X)$ is a square, $p^2$ divides $N$,
and this implies that $N \in \mathcal M_E$.

\medskip (ii) We have $|{\rm det}(S_X)| =dc^2$, and this implies that $S_X \simeq q_I$ where $I$
is an $O$-ideal of an order $O$ of conductor $c$ of $K$. We have $\ell(X) \leqslant 2$, hence $c \in \mathcal M_E$. 

\begin{prop} Suppose that $[E:{\bf Q}] = 20$ and that $S_3 \not = \varnothing$. Set
$d_E = \underset {p \in S_3} \prod p$ and $K_E = {\bf Q}(\sqrt d_E)$. If $c \in \mathcal M_E$, then there exists a complex projective $K3$ surface with maximal complex multiplication by $E$
with Picard lattice isomorphic to $q_I$ for some projective ideal $I$ of an order of conductor $c$ of $K_E$. 

\end{prop}

\noindent
{\bf Proof.} Let $p = p_i$ be a divisor of $c$, and let $P$ be a prime $O_E$-ideal with $p = {\rm N}_{E/{\bf Q}}(P)$ if $p \in S_1$ or $S_3$, 
 and $p^2 = {\rm N}_{E/{\bf Q}}(P)$ if $p \in S_2$. Set $e_P = n_i$, and $G_1 = \underset{P} \bigoplus \ O_E/P^{e_P}$. 
 For all $p_i \in S_3$, let $P_i$ be a prime $O_E$-ideal above $p_i$, and set  $G_2 = \underset{p_i \in S_3} \bigoplus O_E/P_i$; let
 $G = G_1 \oplus G_2$. By Theorem \ref{disc} there
 exists an $O_E$-lattice $T$ of signature $(2,18)$ and discriminant group $G$. The lattice $T$ embeds primitively into
 the $K3$-lattice $\Lambda$, hence by Proposition \ref{existence} there exists a complex projective $K3$ surface $X$ with transcendental
 lattice $T$ and maximal complex multiplication by $O_E$. The orthogonal complement of $T$ in $\Lambda$ is an indefinite binary quadratic
 form of determinant $-d_E c^2$, hence it is of the form $q_I$ for some projective ideal $I$ of an order of conductor $c$ of $K_E$. 

\begin{example} Let $E = {\bf Q}(\zeta_{25})$, and let $P$ be the unique ramified prime $O_E$-ideal. We have $f_P = 1$,
hence with the above notation we have $d_E = 5$ and $K_E = {\bf Q}(\sqrt 5)$.

\medskip
$\bullet$ $p \in S_1 \iff p \equiv 1 \ {\rm (mod \ 25)}$; 

\medskip
$\bullet$ $p \in S_2 \iff p^2 \equiv 1 \ {\rm (mod \ 25)}$.

\medskip
$\bullet$ $S_3 = \{5 \}$.

\medskip Proposition \ref{Picard ramified} implies that if $X$ is a $K3$ surface with maximal complex multiplication by $E$, then
$S_X \simeq q_I$, where $I$ is a projective ideal of an order $O$ of $K_E$. Moreover, if $c$ is the conductor of $O$, then 
we have $c = 5^{2r} N$ where $r \geqslant 0$ is an integer, and if $p$ is a prime divisor of $N$, then $p \equiv \pm 1 \ {\rm (mod \ 25)}$.

\medskip
For $N = 1$, we recover one of  Vorontsov's $K3$ surfaces, see \cite{V}, \cite{LSY}.

\end{example}

\section{$K3$ surfaces with maximal complex multiplication by cyclotomic fields}\label{cyclotomic section}

We keep the notation of the previous sections, and suppose that $E$ is a cyclotomic field. We consider $E$ embedded in ${\bf C}$,
with $E = {\bf Q}(\zeta_m)$, where $m \geqslant 3$ is an integer and $\zeta_m$ is a primitive $m$-th root of unity. As in the previous
sections, the degree of $E$ is denoted by $2n$; note that $2n = \varphi(m)$, and that by hypothesis $2n \leqslant 20$. 

\medskip Recall that if  $X$ is a $K3$ surface with complex multiplication by $O_E$, we denote by $G_X$ the discriminant $O_E$-module $T_X^{\sharp}/T_X$, and that the
minimal number of generators (as an abelian group) of $G_X$ is denoted by $\ell(X)$; it is called the length of $X$. 

\medskip
We start by observing that $K3$ surfaces with maximal complex multiplication by $E$ of length $\leqslant 20-2n$
 are determined by their discriminant modules.

\begin{prop}\label{cyclotomic Nikulin} Let $X$ and $Y$ be two $K3$ surfaces with maximal complex multiplication by $E$ of length 
$\leqslant 20-2n$.
Then $X$ and $Y$ are isomorphic
if and only if the discriminant $O_E$-modules of $T_X$ and $T_{Y}$ are isomorphic.

\end{prop}

\noindent
{\bf Proof.} 
We have $h_E = 1$, since $2n \leqslant 20$ (see for instance \cite{W}, Tables, \S 3
and \S 4), therefore the proposition follows from Proposition \ref{Nikulin}.

\section{Discriminant forms of Craig-like lattices}\label{Craig section}

Let $p$ be a prime number, $p \not = 2$. Craig's lattices are positive definite lattices associated to the cyclotomic field ${\bf Q}(\zeta_p)$; see
for instance \cite{BB}, \S 4 or \cite{CS}, \S 5.4. In this section and the next one, we define (definite and indefinite) analogs of these
lattices.

\medskip
Let $r \geqslant 1$ be an integer, and set $E = {\bf Q}(\zeta_{p^r})$; we have $[E:{\bf Q}] = 2n = (p-1)p^{r-1}$. Let $P$ be the unique prime ideal
of $O_E$ above $p$, let us write $D_E = P^{\delta}$, and let $d$ be an integer such that $\delta = 1 - 2d$. 

\medskip
We consider the quadratic space $(E,q)$ with $q : E \times E \to {\bf Q}$ defined by $q(x,y) = {\rm Tr}_{E/{\bf Q}}(x \overline y)$. For
all integers $k$ such that $k \geqslant d$, let us denote by $C_k$ the lattice of $(E,q)$ given by $C_k = P^k$.

\begin{notation} If $b \not = 0$ is an integer, we denote by $({\bf Z}/p{\bf Z},{b \over p}xy)$ the symmetric bilinear form $${\bf Z}/p{\bf Z} \times {\bf Z}/p{\bf Z} \to {\bf Q}/{\bf Z}$$
sending $(x,y)$ to ${b \over p}xy$.

\end{notation}

Let $W({\bf Q}/{\bf Z})$ be the Witt group of symmetric bilinear forms on finite abelian groups (see for instance \cite{Sch}, Chapter 5, \S 1),
and let $[G,q]$ be the Witt class of $(G,q)$ in $W({\bf Q}/{\bf Z})$.

\begin{theo}\label{Craig discriminant} Let $L = C_k$ with $k$ as above, set $a = \delta + 2k$, and let $e \in \{\pm 1\}$ be such
that $p^{r-1} \ \equiv e\ {\rm (mod \ 4)}.$  We have $G_L \simeq O_E/P^a$, and $$[(G_L,q_L)] = [({\bf Z}/p{\bf Z},{-e \over p}xy)]$$

in $W({\bf Q}/{\bf Z})$.

\end{theo}

\begin{lemma}\label{independent} The Witt class of $(G_{C_k},q_{C_k})$ is independent of $k$.

\end{lemma}

\noindent
{\bf Proof.} Let $k,\ell$ be such that $d \leqslant k \leqslant \ell$. We have $C_{\ell} \subset C_{k}$, hence $C_k^{\sharp} \subset C_{\ell}^{\sharp}$;
therefore $C_k/C_{\ell}$ is totally isotropic in $C_{\ell}^{\sharp}/C_{\ell}$, and $(C_k/C_{\ell})^{\perp} = C_k^{\sharp}/C_{\ell}$.
By \cite{Sch}, Lemma 5.1.3 this implies that the Witt classes of $(G_{C_k},q_{C_k})$ and of $(G_{C_{\ell}},q_{C_{\ell}})$
are equal.

\begin{lemma}\label{dual} We have $C_k^{\sharp} = C_{-\delta -k}$ and $C_{k + 2n} = p C_k$.

\end{lemma}

\noindent
{\bf Proof.} Indeed,  $C_k^{\sharp} = D_E^{-1} P^{-k} = C_{-\delta - k}$, and $C_{k + 2n} = P^k P^{2n} = pC_k$. 

\begin{lemma}\label{p} Suppose that $r = 1$. Then the lattice $C_d$ is isomorphic to the root lattice $A_{p-1}$, and
$C_{-d + 1}$ is isomorphic to $p A_{p-1}^{\sharp}$. 

\end{lemma}

\noindent
{\bf Proof.} We have $C_d \simeq A_{p-1}$ by \cite{E}, Lemma 5.4. The second assertion follows from this, and the previous lemma. 

\begin{prop}\label{p bis} Assume that $r = 1$, and let $L = C_d$. Then $(G_L,q_L)$ is isomorphic to $({\bf Z}/p{\bf Z},{-1\over p}xy)$. 

\end{prop}

\noindent
{\bf First proof of Proposition \ref{p bis}.} It is clear that $G_L \simeq {\bf Z}/p{\bf Z}$. To show that $q_L(x,y) = {-1\over p}xy$, 
apply \cite{BT}, Proposition 6.3 with $\lambda = 1$. Note that (with the notation of \cite{BT}) we have $\pi_E^{-2d+2} = \pi_E^{p-1} = p$,
and that $\overline \pi_E = - \pi_E$. This implies that the invariant unit $u$ of Proposition 6.3 is 
${\rm Tr}_{E/{\bf Q}}({1 \over p})$ if $d$ is even, and ${\rm Tr}_{E/{\bf Q}}({-1 \over p})$ if $d$ is odd; hence 
$u = {{p-1}\over 2}$ if  $p \equiv 3 \ {\rm (mod \ 4)}$
and
$u = {{1-p}\over 2}$ if  $p \equiv 1 \ {\rm (mod \ 4)}$. Suppose first that $p \equiv 3 \ {\rm (mod \ 4)}$. Then 
 $q_L(x,y) = {-1\over p}xy$, as claimed. If $p \equiv 1 \ {\rm (mod \ 4)}$, then $-1$ is a square ${\rm (mod \ 4)}$, hence 
 $q_L(x,y) = {-1\over p}xy$ in this case as well. 
 
 \medskip
 \noindent
 {\bf Second proof of Proposition \ref{p bis}.} By Lemma \ref{p}, the lattice $L = C_d$ is isomorphic to the root lattice $A_{p-1}$, and
 hence $(G_L,q_L)$ is isomorphic to $({\bf Z}/p{\bf Z},{-1\over p}xy)$; see for instance  \cite{Mc2}, Proposition 3.5. 
 
 \begin{remark} Assume that $r = 1$. With the notation of \cite{BB}, we have $C_k = A^{\ell}_{p-1}$, where $\ell = k + {{p-1} \over 2}$. 
 
 \end{remark}
 
 \noindent
 {\bf Proof of Theorem \ref{Craig discriminant}.} 
 We have $L \simeq P^k$ and $L^{\sharp} \simeq P^{-k} D_E^{-1} = P^{-k} P^{-\delta}$, hence $G_L \simeq O_E/P^{\delta + 2k} = O_E/P^a$. 
 To prove that $[(G_L,q_L)] = [({\bf Z}/p{\bf Z},{-e \over p}xy)]$
in $W({\bf Q}/{\bf Z})$.
  we may assume that $L = O_E$ (cf. Lemma \ref{independent}).
  By \cite{B 06}, Proposition 9.1, the lattice $(O_E,q)$ is isomorphic to the orthogonal sum of 
 $p^{r-1}$ copies of $p^r A_{p-1}^{\sharp}$. Set $M = p^r A_{p-1}^{\sharp}$ and $T = ( {\bf Z}/p{\bf Z},{{-1 \over p}}xy)$. The Witt class of $(G_M,q_M)$
 in $W({\bf Q}/{\bf Z})$ is equal to $T$; indeed, Lemma \ref{p} implies that $p A_{p-1}^{\sharp}$
 is isomorphic to $C_{-d + 1}$  for $r = 1$, and hence by Lemma \ref{dual}  the lattice $p^r A_{p-1}^{\sharp}$ is also of the form $C_{\ell}$ for
 some $\ell$ and $r = 1$. By Proposition \ref{p bis} and Lemma \ref{independent}, this implies that $(G_M,q_M)$ and $T$ are
 in the same class in $W({\bf Q}/{\bf Z})$; hence
 $(G_L,q_L)$ is Witt equivalent to the orthogonal sum of $p^{r-1}$ copies of $T$. In $W({\bf Q}/{\bf Z})$, the orthogonal
 sum of $4$ copies of $T$ is always $0$, the sum of two copies of $T$ is $0$ if and only if $p \equiv 1 \ {\rm (mod \ 4)}$.
 This implies that $(G_L,q_L)$ is Witt equivalent to $T = ( {\bf Z}/p{\bf Z},{{-e \over p}}xy)$, 
 as claimed. 
 
 \begin{example} Let $p = 3$ and $r = 2$; then $d = -4$, and the lattice $C_d$ is isomorphic to the root lattice $E_6$
 (see \cite{B 99}, \S 3). 
 
 \end{example}
 
 \section{Indefinite Craig-like lattices}\label{indefinite Craig section}
 
 We keep the notation of the previous section. Let $F$ be the maximal totally real subfield of $E$, and let
 $\sigma_0 : F \to {\bf R}$ be a real embedding of $F$. In this section, we assume that {\it there exists a unit
 $u \in O_F^{\times}$ such that 
 $\sigma_0(u) > 0$ and that $\sigma(u) < 0$ for all the other real embeddings $\sigma$ of $F$}. 
 
\medskip
Let $h^-$ be 
 the relative class number of $E$ (i.e. the class number of $E$ divided by the class number of $F$). 
 If $h^-$ is odd, then there exists a unit $u$ as above, and its image in $O_F^{\times}/{\rm N}_{E/F}(O_E^{\times})$ is unique (see for instance
 \cite {B 84}, Lemma 3.2).

 \medskip
 We consider the quadratic space $(E,q_0)$, where $q_0 : E \times E \to {\bf Q}$ is given by
 $q_0(x,y) = {\rm Tr}_{E/{\bf Q}} (u x \overline y)$; the signature of $(E,q_0)$ is $(2,2n-2)$. 
 
 \medskip
 If $k$ is an integer with $k \geqslant d$, we denote by $L_k$ the lattice of $(E,q_0)$ given by $L_k = P^k$. 
 
 \begin{lemma}\label{independent bis} The Witt class of the discriminant form of $L_k$ in $W({\bf Q}/{\bf Z})$ 
 is independent of $k$.
 
 \end{lemma}
 
 \noindent
 {\bf Proof.} This follows from the same argument as Lemma \ref{independent}. 
 
 \medskip
 Let $\overline \epsilon \in {\bf F}_p^{\times}/{\bf F}_p^{\times 2}$ be the unique non-trivial element, and let $\epsilon \in {\bf Z}$
 be such that the image of $\epsilon$ in ${\bf F}_p^{\times}/{\bf F}_p^{\times 2}$ is equal to $\overline \epsilon$. 
 
 \begin{theo}\label{indefinite Craig discriminant} Let $L = L_k$ be as above, set $a = \delta + 2k$ and let 
 $e \in \{\pm 1\}$ be such
that $p^{r-1} \ \equiv e\ {\rm (mod \ 4)}.$  We have $G_L \simeq O_E/P^a$, and $$[(G_L,q_L)] = [({\bf Z}/p{\bf Z},{\epsilon e \over p}xy)]$$

in $W({\bf Q}/{\bf Z})$.

 \end{theo}

 \noindent
 {\bf Proof.} We have $G_L \simeq O_E/P^a$ as in Theorem \ref{Craig discriminant}. Let us show
 that $[(G_L,q_L)] = [({\bf Z}/p{\bf Z},{\epsilon e \over p}xy)]$
in $W({\bf Q}/{\bf Z})$.
 By Lemma \ref{independent bis} 
 it is enough to consider the case $k = d$, hence we can assume that $G_L = {\bf Z}/p{\bf Z}$. 
 By Theorem \ref{Craig discriminant}, the discriminant form of the positive definite lattice $C_d$ is $({\bf Z}/p{\bf Z},{-e \over p}xy)$; let
 us show that the discriminant form of $L_d$ is $({\bf Z}/p{\bf Z},{\epsilon e \over p}xy)$. 
 
 \medskip
 Recall that $n = [F:{\bf Q}]$. Let $v_p$ be the unique finite place of $F$ above $p$. Let us write $E = F(\sqrt \theta)$, with $\theta \in F^{\times}$. 
 We have $(u,\theta)_v = 1$ for all finite places of $v$ of $F$ with $v \not = v_p$. Suppose that $n$ is odd. Then  $(u,\theta)_v = 1$
 at an even number of infinite places $v$ of $F$, hence the product formula implies that $(u,\theta)_{v_p} = 1$.
 This implies that the discriminant forms of $L_d$ and $C_d$ are isomorphic, hence they are isomorphic
 to $({\bf Z}/p{\bf Z},{-e \over p}xy)$. Note that $n$ is odd if and only if $p \equiv 3 \ {\rm (mod \ 4)}$, and in this case
 $-1$ is not a square ${\rm (mod \ p)}$, hence we can take $\epsilon = -1$. 
 
 \medskip
 Suppose that $n$ is even. Then the above argument shows that $(u,\theta)_{v_p} = -1$. Applying
 \cite{BT}, Proposition 6.6, we see that the discriminant form of $L_d$ is isomorphic to $({\bf Z}/p{\bf Z},{- \epsilon e \over p}xy)$
 in this case. Note that $n$ is even if and only if $p \equiv 1 \ {\rm (mod \ 4)}$, and $-1$ is a square 
 ${\rm (mod \ 4)}$ in this case. This implies that the discriminant form
 of $L_d$ is isomorphic to $({\bf Z}/p{\bf Z},{\epsilon e \over p}xy)$ in this case as well.
 
 \begin{remark}
 The lattices $L_k$ and $C_k$ are defined for $k \geqslant d$, and their determinant is $p^a$, with $a = \delta + 2k$. The condition
 $k \geqslant d$ is equivalent to $a \geqslant 1$. This motivates the following notation :
 
 \end{remark}
 
 \begin{notation}\label{lambda}  If $a \geqslant 1$ is an integer, set $\Lambda_a = L_k$  with $k = {{a - \delta}\over 2}$.
 
 \end{notation}
 
 \begin{example}\label{first example} Let $E = {\bf Q}(\zeta_p)$, i.e. $r = 1$. Then $e = 1$, hence the discriminant form of
 $\Lambda_1$ is $({\bf Z}/p{\bf Z},{\epsilon  \over p}xy)$. If $p \equiv 3 \ {\rm (mod \ 4)}$, then we can take
 $\epsilon = -1$, and this implies that the discriminant form is $({\bf Z}/p{\bf Z},{-1 \over p}xy)$.
 Suppose that $p \equiv 1 \ {\rm (mod \ 4)}$. If $p \equiv 5 \ {\rm (mod \ 8)}$, then $2$ is not a square 
 ${\rm (mod \ p)}$, therefore we can take $\epsilon = 2$, and the discriminant form is then isomorphic
 to $({\bf Z}/p{\bf Z},{2  \over p}xy)$. This covers all the prime numbers $p$ needed for the applications to
 $K3$ surfaces, except for $p = 17$. In this case, we can take $\epsilon = 3$, hence the discriminant
 form is isomorphic to $({\bf Z}/17{\bf Z},{3\over 17} xy)$. 
 
 \end{example}
 
 \begin{example}\label{second example} Let $E = {\bf Q}(\zeta_m)$, with $m = 9, 25$ or $27$. We have
 $e = -1$ if $m = 9$ and $e = 1$ for $m = 25$ or $27$. We can take $\epsilon = -1$ for $m = 9$ or $27$, and
 $\epsilon = 2$ if $m = 25$ (cf. Example \ref{first example}). Hence the discriminant form of $\Lambda_1$ is
 
 \medskip
 $({\bf Z}/p{\bf Z},{1 \over 3}xy)$ if $m = 9$;
 $({\bf Z}/p{\bf Z},{-1 \over 3}xy)$ if $m = 27$;
 $({\bf Z}/p{\bf Z},{2 \over 5}xy)$ if $m = 25$.

 \end{example}
 
 \medskip
 
 If $L$ is a lattice and $n$ is an integer, we denote by $L(n)$ the lattice with values multiplied by $n$. 
 
 \begin{notation}\label{delta}  If $a \geqslant 1$ is an integer, set $\Delta_a = C_k(-1)$  with $k = {{a - \delta}\over 2}$.
 
 \end{notation}
 
 \begin{lemma}\label{minus discriminant} Suppose that $r = 1$, and let $a \geqslant 1$ be an integer. Set $L = \Delta_a$.  Then we have
 $$[(G_L,q_L)] = [({\bf Z}/p{\bf Z},{1 \over p}xy)]$$
in $W({\bf Q}/{\bf Z})$.

 \end{lemma}
 
 \noindent
 {\bf Proof.} This follows from Theorem \ref{Craig discriminant}, noting that since $r = 1$, we have $e = 1$, and that $\Delta_a$ is negative
 definite.

 \section{Twisted lattices}
 
 We keep the notation of the previous section, and assume that $r = 1$, hence $E = {\bf Q}(\zeta_p)$, and $F = {\bf Q}(\zeta_p + \zeta_p^{-1})$. 
 We fix a real embedding $\sigma_0 : F \to {\bf R}$, and we assume that  there exists a unit
 $u \in O_F^{\times}$ such that 
 $\sigma_0(u) > 0$ and that $\sigma(u) < 0$ for all the other real embeddings $\sigma$ of $F$. For all odd integers $a \geq 1$, we
 define the lattices $\Lambda_a$  and $\Delta_a$ as in the previous section. We denote by $P$ the unique ramified ideal of $O_E$. 
 
 \begin{defn} Let $(L,q)$ be a negative definite even lattice. A {\it root} of $L$ is an element $x \in L$ such that $q(x,x) = -2$. 
 We say that $L$ is a root lattice if it has an integral basis of roots. 
 
 \end{defn}

 \begin{prop}\label{root} Let $J \subset O_E$ be an ideal, and let $a \geqslant 1$ be an odd integer; let
 $L$ be a twist of $\Delta_a$ by $J$. If $a > 1$ or if $J \not = O_E$, then the lattice $L$ does not contain any roots.
 
 \end{prop}
 
 \noindent
 {\bf Proof.} The lattice $L$ has an isometry $t : L \to L$ with characteristic polynomial $\Phi_p$. By \cite{BM}, Appendix, this shows that $L$ does not contain
 any root sublattice if $a > 1$ or if $J \not = O_E$.

 \begin{prop}\label{isometry} Let $p$ be a prime number with $p \equiv 3 \ {\rm (mod \ 4)}$, and let $\zeta$ be a primitive $p$-th root of
 unity. Let $a \geqslant 1$ be an integer, and let $J \subset O_E$ be an ideal prime to $P$.
 Let $\Lambda$ be a twist of $\Lambda_a$ by $J$, and let $\Delta$ be a twist of  
$ \Delta_a$ by $J$.
 
 \medskip Let 
 $t_{\Lambda} : \Lambda \to \Lambda$ and $t_{\Delta} : \Delta \to \Delta$ be the
 isometries induced by multiplication by $\zeta$. 
 Then there exists an even, unimodular lattice $L$ containing $\Lambda \oplus \Delta$
 as a sublattice of finite index, and an isometry $t : L \to L$ such that 
$t|\Lambda = t_{\Lambda}$, $t|\Delta = t_{\Delta}$. 
 
 \end{prop}
 
 \noindent
 {\bf Proof.} The proof uses gluing of lattices and isometries, as in McMullen's papers \cite{Mc2}, \S 2 or \cite{Mc3}, \S 4. We check that the conditions of \cite{Mc2}, page 5
 (gluing of a pair of lattices, extending isometries) hold for $\Lambda$, $t_{\Lambda}$ and $\Delta$, $t_{\Delta}$. If suffices to check these conditions for all
 prime numbers $q$ dividing the orders of the discriminant groups (glue groups) of $\Lambda$ and $\Delta$, because of the primary decomposition of 
 these groups (see \cite{Mc2}, page 4). 
 
 \medskip We start with the $p$-primary components. The discriminant modules of $\Lambda_a$ and $\Delta_a$ are isomorphic, and
 their discriminant forms have opposite signs, since their Witt classes are $[({\bf Z}/p{\bf Z},{-1 \over p}xy)]$, respectively $[({\bf Z}/p{\bf Z},{1 \over p}xy)]$;
 hence McMullen's conditions hold at $p$.
 
 \medskip 
 Let $q$ be a prime number with $q \not = p$. If $J$ is squarefree, then the discriminant groups are vector spaces over ${\bf F}_q$, hence the conditions of \cite{Mc2}, \S 3 are fulfilled; if $J$ is not squarefree, then Milnor's argument in \cite{M}, \S 3, Theorem 3.4, shows that it is
 enough to consider the case where the discriminant groups are  vector spaces over ${\bf F}_q$.

 \medskip Hence McMullen's conditions hold at every prime number, and therefore
 there exists an even, unimodular lattice $L$ containing $\Lambda \oplus \Delta$ with finite index, and an isometry $t : L \to L$
 such that $t|\Lambda = t_{\Lambda}$, $t|\Delta = t_{\Delta}$. 

\begin{remark} Let $a = 1$ and $J = O_E$. The previous proposition holds with $t_{\Delta} = id$;
we obtain an even, unimodular lattice $L$ and an isometry $t : L \to L$ such that $t|\Lambda$ is the multiplication by
$\zeta$ and $t|\Delta$ is the identity. 


\end{remark}

\section{A family of $K3$ surfaces}\label{family}

Let $p$ be a prime number, $3 \leqslant p \leqslant 11$. We keep the notation of the previous sections; in particular, $\epsilon \in {\bf Z}$
 is such that the image of $\epsilon$ in ${\bf F}_p^{\times}/{\bf F}_p^{\times 2}$ is the unique non-trivial element.  Set $E = {\bf Q}(\zeta_p)$,
 and let $P$ be the unique ramified prime ideal of $O_E$.

\medskip
If $X$ is a $K3$ surface, we denote by $T_X$ its transcendental lattice and by $S_X$ its Picard lattice.

 \begin{theo}\label{K3 family} Let $a \geqslant 1$ be an  odd integer.
 There exists a unique {\rm (up to isomorphism)} complex projective
 $K3$ surface $X_a = X_a(p)$ with maximal complex multiplication by $E$ such that the following equivalent conditions hold~:
 
 \medskip
 {\rm (i)}  The discriminant module of  $T_{X_a}$ is isomorphic to $P/P^a$, and the Witt class of its discriminant form $(G_{X_a},q_{X_a})$ is 
 $[({\bf Z}/p{\bf Z},{\epsilon  \over p}xy)]$.

  \medskip
 {\rm (ii)} $G_{X_a} \simeq P/P^a$.
 
  \medskip
 {\rm (iii)} ${\rm det}(S_{X_a}) = {\rm det}(T_{X_a}) = p^a$. 
 
 \medskip Moreover, the surfaces $X_a$ are isogeneous for all $a \geqslant 1$. 
 
 \end{theo}
 
 \medskip
 
 Recall that for all integers $a \geqslant 1$, the $O_E$-lattice $\Lambda_a$ is defined in \S \ref{indefinite Craig section}, see
 notation \ref{lambda}. 
 
 \begin{theo}\label{K3 lemma} Let $a \geqslant 1$ be an  odd integer.
 There exists a unique {\rm (up to isomorphism)} complex projective
 $K3$ surface $X_a$ with maximal complex multiplication by $E$ such that the transcendental lattice $T_{X_a}$ is Hodge isomorphic
 to the $O_E$-lattice  $\Lambda_a$. 
 
 \end{theo}
 
 \noindent
 {\bf Proof.} Since $p \leqslant 11$, we have $\ell(G_{\Lambda_a}) \leqslant 10$,  hence
 $\Lambda_a$ embeds uniquely  into $\Lambda$. By Proposition \ref{existence}, there exists a complex projective $K3$ surface $X_a$ having CM by $E$ such that the
 $O_E$-lattice $T_{X_a}$ is Hodge isometric to $\Lambda_a$; the surface $X_a$ is unique up to isomorphism (cf. Proposition 
 \ref{iso surfaces}). By construction, we have $\Lambda_a \otimes _ {\bf Z} {\bf Q} \simeq \Lambda_b \otimes _ {\bf Z} {\bf Q}$ for
 all $a, b \geqslant 1$, hence these surfaces are all isogeneous (see \cite{Mu}, \cite{N 87}, \cite{Bu}). 
 
 \medskip
 \noindent
 {\bf Proof of Theorem \ref {K3 family}.} The existence and the uniqueness of the $K3$ surfaces $X_a$ with property (i) results from Theorem \ref{K3 lemma}
and Theorem \ref{indefinite Craig discriminant}; Theorem \ref{K3 lemma} also implies that the surfaces $X_a$ are all isogeneous. 
It is clear that (i) $\implies$ (ii) $\implies$ (iii), and by Proposition \ref{Nikulin} we have (ii) $\implies$ (i). It remains to show that (iii) implies (ii). 
The only prime ideal of $O_E$ above $p$ is the unique ramified ideal $P$, and we have ${\rm N}(P) = p$. Therefore ${\rm det}(T_X) = p^a$
implies that $G_{T_X} \simeq O_E/P^a$ as $O_E$-modules, and therefore (ii) holds. 

\begin{example}\label{Vor 1} For $a = 1$, we recover Vorontsov's examples of $K3$ surfaces (see \cite{V}, \cite{Ko}) arising in
connection with automorphisms acting trivially on the Picard lattice; these surfaces are elliptic with a section and defining
Weierstrass equations over $\bf Q$ are given in  \cite{Ko}.  The
discriminant forms of these surfaces are computed using an elliptic fibration in  \cite{LSY}, Table 5.

\end{example}

 \section{More $K3$ surfaces}\label{more}
 
 Let $E = {\bf Q}(\zeta_m)$ with $m = 9, 13, 17, 19, 25$ or $27$, and let us denote by $P$ the unique ramified ideal of $O_E$. 
 
 \begin{theo}\label{K3 bis} Let $a \geqslant 1$ be an  odd integer, and suppose that $a \leqslant 7$ if $m = 13$, $a \leqslant 5$ if $m = 17$,
 $a =  1$ if $m = 19, 25$ or $27$.
 
 \medskip
 Then there exists a unique {\rm (up to isomorphism)} complex projective
 $K3$ surface $X_a(m)$ with maximal complex multiplication by $E$ such that the following equivalent conditions hold :
 
 \medskip
 {\rm (i)} The discriminant module of  $T_{X_a}$ is isomorphic to $P/P^a$, and the Witt class of its discriminant form $(G_{X_a},q_{X_a})$ is 
 $[({\bf Z}/p{\bf Z},{\epsilon  \over p}xy)]$ if $m = p = 13, 17$ or $19$; $[({\bf Z}/3{\bf Z},{1  \over 3}xy)]$ if $m = 9$; 
 $[({\bf Z}/3{\bf Z},{-1  \over 3}xy)]$ if $m = 27$ and $[({\bf Z}/5{\bf Z},{2  \over 5}xy)]$ if $m = 25$.

  \medskip
 {\rm (ii)} $G_{X_a(m)} \simeq  P/P^a$.
 
 \medskip
 {\rm (iii)} ${\rm det}(S_{X_a}(m)) = {\rm det}(T_{X_a}(m)) = p^a$ if $m$ is a power of $p$. 
 
 \medskip Moreover, the surfaces $X_a(m)$ are isogeneous for all $a \geqslant 1$. 
 
 \end{theo}
 
 \noindent
 {\bf Proof.} Suppose first that $m \not = 25$. Then the proof goes along the same lines as the proof of Theorem \ref{K3 family}; the conditions on $a$ ensure that 
  $\Lambda_a$ embeds uniquely in $\Lambda$ (see \cite{N}, Corollary 1.12.3 and Theorem 1.14.4). By Proposition \ref{existence} and
  Proposition \ref{iso surfaces}), there exists a 
  unique (up to isomorphism) complex projective $K3$ surface $X_a$ having CM by $E$ such that the 
 $O_E$-lattice $T_{X_a}$ is Hodge isometric to $\Lambda_a$
  of \S \ref{indefinite Craig section}, cf. Notation \ref{lambda}. 
  
  \medskip If $m = 25$, then  $\Lambda_a$ embeds primitively in $\Lambda$ (see \cite{N}, Corollary 1.12.3), hence there exists a complex projective $K3$ surface $X_a$ having CM by $E$ such that the 
 $O_E$-lattice $T_{X_a}$ is Hodge isometric to $\Lambda_a$. The surface $X_a$ has an automorphism that is the identity on the Picard lattice
 and induces complex multiplication on the transcendental lattice; this can be checked by noting that the action of the complex multiplication on
 $T_{X_a}$ induces multiplication by -1 on $G_{X_a}$. 
  
  \medskip 
 The discriminant form of $X_a(m)$ can be deduced from  Theorem \ref{indefinite Craig discriminant}, 
 see also Example \ref{second example}; the proof is completed as in the proof of Theorem \ref{K3 family}.

 \begin{example}\label{Vor 2} As in Example \ref{Vor 1}, for $a = 1$ we recover Vorontsov's examples of $K3$ surfaces (see \cite{V}, \cite{Ko}).
 These surfaces are elliptic with a section, except if  $m = 25$; see \cite {Ko} for defining equations over $\bf Q$. 
   
 \end{example}
 
 \section{Automorphisms}\label{automorphisms}
 
 Let $p = 3, 7$ or $11$ and let $a \geqslant 1$ be an odd integer; let $X_a = X_a(p)$ be the $K3$ surface
 of \S \ref{family}. Let $\zeta$ be a primitive $p$-th root of unity, and let $E = {\bf Q}(\zeta)$, considered as a subfield of $\bf C$. 
 
 \begin{defn} Let $X$ be a complex projective $K3$ surface with  complex multiplication by $E$. Let $T : X \to X$ be an automorphism, and 
 let $t : T_X \to T_X$ be the isometry induced by $T$. 
 We say that $t$ {\it induces the complex multiplication by} $E$ if $t(\sigma) = \zeta \sigma$ where $\sigma$ is a non-zero $2$-form in $ T_X \otimes_{\bf Q} {\bf C}$. 
 
 \end{defn}
 
 \begin{theo}\label{no J} {\rm (i)} For all $a \geqslant 1$, the $K3$ surface $X_a$ has an automorphism of order $p$ inducing the complex multiplication by $E$. 
 
 \medskip
 {\rm (ii)} For all $a \geqslant 1$, the surface $X_a$ is elliptic with a section.
 
 \medskip
 {\rm (iii)} If $a > 1$, then the Mordell-Weil lattice of $X_a$ is isomorphic to $\Delta_a$. 
 
 \end{theo}

 \noindent {\bf Proof.} If $a = 1$, then this is well-known : the surfaces $X_1$ are isomorphic to  Vorontsov's $K3$ surfaces, 
 see Example \ref{Vor 1}. 
 
 \medskip
 Suppose that $a > 1$. 
 By Proposition \ref{isometry} there exists an even, unimodular lattice $L$ containing $\Lambda_a \oplus \Delta_a$
 as a sublattice of finite index, and an isometry $t_L: L \to L$ such that 
$t_L|\Lambda_a = t_{\Lambda}$, $t_L|\Delta_a = t_{\Delta}$. 

\medskip
 Set $r = 24-2p$, and let $M$ be an even, unimodular lattice of signature $(1,r-1)$;
 such a lattice exists (and is unique up to isomorphism) since $r-2$ is divisible by $8$. Let $t_M : M \to M$ be the identity. 
 
 \medskip
 Set $N = M \oplus L$, and let $t : N \to N$ be such that $t|L = t_L$, $t|M = t_M = id$. The lattice $N$ is even, unimodular of
 signature $(3,19)$. Set $S = M \oplus \Delta_a$ and $T = \Lambda_a$. 
 Since $a > 1$, the lattice $\Delta_a$ does not contain any roots (cf. Lemma \ref{root}). 
 The isometry $t$ is the
 identity on $N$, hence $t_S$ satisfies the conditions of McMullen in \cite{Mc3}, \S 6. Therefore by \cite{Mc3}, Theorem 6.1 there exists a complex projective $K3$ surface $X$ with $S_X \simeq S$, $T_X \simeq T$, and
 an automorphism $T : X \to X$ such that $T^* = t$. This $K3$ surface is isomorphic to $X_a$. This shows (i).
 
 \medskip
 The lattice $M$ is even, unimodular and of signature $(1,r)$, hence it has an orthogonal factor isomorphic to the 2-dimensional
 hyperbolic lattice $U$. This implies that $X_a$ is elliptic with a section (see \cite{MWL},  Theorem 11.24)
hence (ii) holds. Note that the orthogonal complement of $U$ in $M$ is a (negative) root lattice.

 \medskip
 If $a > 1$, then $\Delta_a$ has no roots, hence the trivial lattice of the fibration is isomorphic to $M$, and the Mordell-Weil lattice
 to $\Delta_a$; this implies (iii). 
 
 \section{Twisted $K3$ surfaces}\label{twisted K3 section}
 
 The aim of this section is to extend the results of \S \ref{family} and \S \ref{automorphisms} to certain twisted $K3$ surfaces. 
 Let $p$ be a prime number, $3 \leqslant p \leqslant 11$, and set $E = {\bf Q}(\zeta_p)$. We keep the notation of the
 previous sections; in particular, $P$ is the unique ramified prime ideal of $E$. 
 
 \medskip
 If $a \geqslant 1$ is an odd integer, the lattices $\Lambda_a$ and $\Delta_a$ are defined in \S \ref{family}. Let $J \subset O_E$ 
 be an $O_E$-ideal prime to $P$ such that $\overline J = J$. Since $h_E = 1$, the twist of an $O_E$-lattice by $J$ is uniquely defined (up to isomorphism; see Proposition \ref{twisting proposition}). We
 denote by $\Lambda_{a,J}$ and $\Delta_{a,J}$ the twists of $\Lambda_a$ and $\Delta_a$ by $J$. 
 
 \begin{theo}\label{twisted K3 lemma} Let $a \geqslant 1$ be an  odd integer, and let $J \subset O_E$ be an $O_E$-ideal prime to $P$
such that $\overline J = J$. Then

\medskip
 {\rm{(i)}} There exists a unique {\rm (up to isomorphism)} complex projective
 $K3$ surface $X_{a,J}$ with maximal complex multiplication by $E$ such that the transcendental lattice 
 $T_{X_{a,J}}$ is Hodge isomorphic
 to the $O_E$-lattice  $\Lambda_{a,J}$. 
 
 \medskip
 Suppose that $p = 3$, $7$ or $11$. 
 
 \medskip
 {\rm (ii)} The $K3$ surface $X_{a,J}$ has an automorphism of order $p$ inducing the complex multiplication by $E$.
 
 \medskip
 {\rm (iii)} The surface $X_{a,J}$ is elliptic with a section. 
 
 \medskip
 {\rm (iv)} If $a > 1$ or $J \not = O_K$, then the Mordell-Weil lattice of $X_{a,J}$ is isomorphic to $\Delta_{a,J}$. 
 
 \end{theo}
 
 \noindent
 {\bf Proof.} (i)
 We have $[E:{\bf Q}] = p-1$ and $p \leqslant 11$, hence the lattice
 $\Lambda_{a,J}$ embeds uniquely into $\Lambda$. Proposition \ref{existence} implies that  there exists a complex projective $K3$ surface $X_{a,J}$ having CM by $E$ such that the
 $O_E$-lattice $T_{X_{a,J}}$ is Hodge isometric to $\Lambda_{a,J}$; the surface $X_{a,J}$ is unique up to isomorphism (cf. Proposition 
 \ref{iso surfaces}). 
 
 \medskip
 (ii) and (iii)  If $a = 1$, then this follows from Theorem \ref{no J},  (i) and (ii). Set $\Lambda = \Lambda_{a,J}$ and $\Delta = \Delta_{a,J}$. 
 By Proposition \ref{isometry} there exists an even, unimodular lattice $L$ containing $\Lambda  \oplus \Delta$
 as a sublattice of finite index, and an isometry $t_L: L \to L$ such that 
$t_L|\Lambda = t_{\Lambda}$, $t_L|\Delta = t_{\Delta}$. 
 Set $r = 24-2p$, and let $M$ be an even, unimodular lattice of signature $(1,r-1)$;
 such a lattice exists (and is unique up to isomorphism) since $r-2$ is divisible by $8$. Let $t_M : M \to M$ be the identity. 
Set $N = M \oplus L$, and let $t : N \to N$ be such that $t|L = t_L$, $t|M = t_M = id$. The lattice $N$ is even, unimodular of
 signature $(3,19)$. Set $S = M \oplus \Delta$ and $T = \Lambda$. 
 
 \medskip
 Suppose that $a > 1$. Then the lattice $\Delta$ does not contain any roots (cf. Lemma \ref{root}). 
 The isometry $t$ is the
 identity on $N$, hence $t_S$ satisfies the conditions of McMullen in \cite{Mc3}, \S 6. Therefore by \cite{Mc3}, Theorem 6.1 there exists a complex projective $K3$ surface $X$ with $S_X \simeq S$, $T_X \simeq T$, and
 an automorphism $T : X \to X$ such that $T^* = t$. This $K3$ surface is isomorphic to $X_{a,J}$, and this implies (ii).
 
 \medskip
 The lattice $M$ has an orthogonal factor isomorphic to the $2$-dimensional hyperbolic lattice $U$, therefore $X_{a,J}$ is elliptic
 with section (see \cite{MWL},  Theorem 11.24), hence (iii) holds. 
 
 \medskip
 (iv) If $a > 1$ or $J \not = O_E$, then by Lemma \ref{root} the lattice $\Delta$ does not contain any roots. Therefore the trivial
 lattice of the fibration is isomorphic to $M$, and the Mordell-Weil lattice is isomorphic to $\Delta$. 
 
 \begin{theo}\label{all surfaces} Let $p = 3, 7$ or $11$, and let $X$ be a $K3$ surface having an 
 automorphism of order $p$ inducing the complex multiplication by $E$. Then there exists an integer $a \geqslant 1$ and
 an ideal $J \subset O_E$ such
 that $X \simeq X_{a,J}$.

 \end{theo}
 
 \noindent
 {\bf Proof.} The $O_E$-module $G_X = T_X^{\sharp}/T_X$ is isomorphic to $O_E/P^a \oplus O_E/J$ for some
 integer $a \geqslant 1$ and some $O_E$-ideal $J$.  The discriminant module of $X_{a,J}$ is also isomorphic to
 $O_E/P^a \oplus O_E/J$, and the length of this abelian group is $\leqslant 10$, hence by Proposition \ref{Nikulin} we have
 $X \simeq X_{a,J}$.

 \section{Moduli spaces}\label{moduli}
 
An automorphism of a $K3$ surface is said to be {\it symplectic} if it induces the identity
 on the transcendental lattice (hence on the symplectic form), and {\it non-symplectic} otherwise. Artebani, Sarti and Taki
 identified the irreducible components of the moduli space of $K3$ surfaces with a non-symplectic automorphism of prime order
 (see \cite{AST}, see also \cite{AS} for $p = 3$ and \cite{OZ}, \cite{ACV} for $p = 11$). 

 
 \medskip
 Let $\Lambda = \Lambda_{3,19}$ be the $K3$-lattice.  If $X$ is a $K3$-surface, we denote by $\omega_X$ a nowhere
 vanishing 
 holomorphic 2-form on $X$. Let $p$ be a prime number, and let $\zeta_p$ be a primitive $p$-th root of unity. Let $\rho : \Lambda \to \Lambda$
 be an isometry of order $p$, and let us denote by $[\rho]$ its conjugacy class in $O(\Lambda)$. A $[\rho]$-polarized 
 $K3$ surface is a pair $(X,t)$ where $X$ is a $K3$ surface and $t$ a non-symplectic automorphism of $X$ of
 order $p$ such that $t^*(\omega_X) = \zeta_p \omega_X$ and $t^* = \Phi \circ \rho \circ \Phi^{-1}$, for some (fixed) isometry
 $\Phi : \Lambda \to H^2(X,{\bf Z})$, which is called a marking; the moduli space of such polarized $K3$ surfaces
 is denoted by  ${\mathcal M}^p$  (see \cite{AST}, \S 9). We say that a point of this moduli space is a {\it CM point} if 
 the $K3$ surface $X$ has complex multiplication.

 \medskip
 Let $p = 3, 7$ or $11$, and let $E = {\bf Q}(\zeta_p)$. We keep the notation 
 of the previous sections : in particular, $P$ denotes the unique ramified ideal of $E$.

 \medskip
  Let $a \geqslant 1$  be an  odd integer and let
 $J \subset O_E$ be an ideal prime to $P$.
 Let $X_{a,J}$ be the $K3$ surface defined in \S \ref{twisted K3 section}.  The following
 corollary is an immediate consequence of Theorem \ref{twisted K3 lemma}, (ii) : 
 
 \begin{coro} The $K3$ surface $X_{a,J}$ has an automorphism of order $p$ inducing the complex conjugation on
 the transcendental lattice $T_{X_{a,T}}$, and determines a $CM$ point on the moduli space ${\mathcal M}^p$.
 
 \end{coro}
 
 Moreover,  Theorem \ref{all surfaces} implies that all $CM$ points with maximal complex multiplication by $E$ arise in this way. 
 
 \section{Fields of definition, class fields and elliptic fibrations}\label{fields}
 
 Piatetski-Shapiro and Shafarevich proved that a $K3$ surface with complex multiplication can be defined over a number field (see \cite{PS},
 Theorem 4). 
 If moreover
 the complex multiplication is {\it maximal}, Valloni obtained more precise results in \cite{V 21}, \cite{V 23}; if
 $E$ is a CM field and $I \subset O_E$ an ideal with $\overline I = I$, he defined a finite abelian extension $F_I(E)$ such
 that every $K3$ surface with CM by $O_E$ with discriminant ideal $I$ can be defined
 over $F_I(E)$.
 
 \medskip
 Let $p = 3,7$ or $11$, let $\zeta_p$ be a primitive $p$-th root of unity and let $E = {\bf Q}(\zeta_p)$. 
 Let $a \geqslant 1$ be an integer, 
 let $J \subset O_E$ be an ideal relatively prime to the unique ramified ideal $P$ of $E$ such that
 $\overline J = J$.
 Let $X_{a,J}$ be the $K3$ surface defined in \S \ref{twisted K3 section}. As we have seen in \S \ref{moduli}, 
 this gives rise to a CM point on the moduli space $\mathcal M^p$; the description of Artebani, Sarti and Taki
 of the moduli space can be used to obtain a field of definition of $X_{a,J}$. This is illustrated
 by the following example, due to Brandhorst and Elkies :
 
 \begin{example}\label{BE} McMullen proved the existence of a $K3$ surface with an automorphism of entropy
 equal to the logarithm of the Lehmer number (see \cite{Mc3}), and raised the question of constructing
 this surface and the automorphism explicitly. This was achieved by Brandhorst and Elkies in \cite{BE}. 
 
 \medskip Set $E = {\bf Q}(\zeta_7)$, $F = {\bf Q}(\zeta_7 + \zeta_7^{-1})$ and let $J$ be one of the $O_E$-ideals above $13$. The $K3$ surface
$S$  constructed in \cite{BE} has a non-symplectic automorphism of order 7, and the construction shows that it is isomorphic to $X_{1,J}$. In \cite{BE}, \S 3 and \S 4, the authors use 
the description of $\mathcal M^7$ by \cite{AST} to obtain an equation for the surface $S$, with coefficients in
a quadratic extension $K$ of the field $F$; set $K = {\bf Q}(w)$, where $w$ is such that
$$w^6 - 2w^5 + 2w^4 - 3w^3 + 2w^2 - 2w + 1 = 0.$$ The field $K$ has discriminant $7^413$, and contains the field $F$ of discriminant $49$. 

\medskip
Let $P$ be the unique ramified prime ideal of $E$. The composite field $KE$ is isomorphic to Valloni's number  field $F_{PJ}(E)$, i.e. $F_{PJ}(E) \simeq {\bf Q}(\zeta_7,w)$. This computation was
done with the help of PARI/GP.
 
 \end{example}

\section{Some equations}\label{equations section}
Let $p = 3,7$ or $11$,
let $\zeta_p$ be a primitive $p$-th root of unity, and set $E = {\bf Q}(\zeta_p)$.  
Let $a$ be an odd integer,
let $J \subset O_E$ be an ideal relatively prime to the unique ramified ideal $P$ of $E$ such that
 $\overline J = J$.
 We denote by $X_{a,J}(p)$  the $K3$ surface defined in \S \ref{twisted K3 section}; if $J = O_E$, then we use
 the notation $X_a(p)$, as in \S \ref{family}. As mentioned in Example \ref{BE}, Brandhorst and
 Elkies gave an explicit equation for the surface $X_{1,J}(7)$, where $J$ is one of the prime $O_E$-ideals above 13 (see 
 \cite{BE}, \S 4); their method
can be used for other choices of $a$ and $J$.

\medskip
I thank Simon Brandhorst for the following examples.

\begin{example} \label{7,2}
Let $p = 7$, $E = {\bf Q}(\zeta_7)$,  set $w = \zeta_7 + \zeta_7^{-1}$, and let $J$ be one of the prime $O_E$-ideals above $2$. An equation of the
surface $X_{1,J}(7)$ is given by 

$$y^2 = x^3 + bx + ct^7 + d,$$ where

\medskip

$b =  (-3403/16)(w^2 + 4w + 4) $

\medskip
$c = 14(w^2+2w+1)$

\medskip
$d = (293419/32)(w^2+ 2w+ 1)$.

\end{example}

\begin{example}\label{2,7}
Let $p = 7$, $E = {\bf Q}(\zeta_7)$, and set $w = \zeta_7 + \zeta_7^{-1}$. An equation of the surface $X_3(7)$ is given by 

$$y^2 = x^3 + bx + ct^7 + d,$$ where

\medskip

$b = (-230578777287775/2 )w^2 + (127961567541885/2)w + 4144846476936445/16$

\medskip
$c = -5842669785012830924 w^2 + 3242437110294043228 w + 13128359838180149367$

\medskip
$d = 151461887453084383247079/32$.

\end{example}

\bigskip
The computations are due to Simon Brandhorst, and they were done by Sage.

\bigskip
\bigskip
Eva Bayer--Fluckiger 

EPFL-FSB-MATH

Station 8

1015 Lausanne, Switzerland

\medskip

eva.bayer@epfl.ch

\end{document}